\documentclass[11pt]{article}
\usepackage{psfig}
\usepackage{latexsym}
\usepackage{amssymb}

\def\X{{\cal X}}

\def\S{{\cal S}}
\def\Var{\mbox{\rm Var}}

\def\beq{\begin{eqnarray}}
\def\eeq{\end{eqnarray}}

\begin{document}

\fontsize{11}{14.5pt}\selectfont

\begin{center}

{\small Technical Report No.\ 0406,
 Department of Statistics, University of Toronto}

\vspace*{0.9in}

{\LARGE\bf Improving Asymptotic Variance of MCMC Estimators:\\[6pt]  
           Non-reversible Chains are Better }
\\[16pt]

{\large Radford M. Neal}\\[3pt]
 Department of Statistics and Department of Computer Science \\
 University of Toronto, Toronto, Ontario, Canada \\
 \texttt{http://www.cs.utoronto.ca/$\sim$radford/} \\
 \texttt{radford@stat.utoronto.ca}\\[10pt]

 15 July 2004
\end{center}

\vspace{8pt} 

\noindent \textbf{Abstract.}  I show how any reversible Markov chain
on a finite state space that is irreducible, and hence suitable for
estimating expectations with respect to its invariant distribution,
can be used to construct a non-reversible Markov chain on a related
state space that can also be used to estimate these expectations, with
asymptotic variance at least as small as that using the reversible
chain (typically smaller).  The non-reversible chain achieves this
improvement by avoiding (to the extent possible) transitions that
backtrack to the state from which the chain just came.  The proof that
this modification cannot increase the asymptotic variance of an MCMC
estimator uses a new technique that can also be used to prove Peskun's
(1973) theorem that modifying a reversible chain to reduce the
probability of staying in the same state cannot increase asymptotic
variance.  A non-reversible chain that avoids backtracking will often
take little or no more computation time per transition than the
original reversible chain, and can sometime produce a large reduction
in asymptotic variance, though for other chains the improvement is
slight.  In addition to being of some practical interest, this
construction demonstrates that non-reversible chains have a
fundamental advantage over reversible chains for MCMC estimation.
Research into better MCMC methods may therefore best be focused on
non-reversible chains.

\section{\hspace*{-7pt}Introduction}\label{sec-intro}\vspace*{-10pt}

Markov chain Monte Carlo (MCMC) is widely used to estimate
expectations of functions with respect to complex, high-dimensional
probability distributions, particularly in Bayesian statistics and
statistical physics (see, for example, Liu 2001).  An MCMC estimator
can be based on any Markov chain that is irreducible and that has the
distribution of interest as its invariant distribution.  However, the
choice of Markov chain will affect the efficiency with which estimates
of expectations with a given accuracy can be obtained.  In this paper,
I show that an MCMC estimator based on a reversible Markov chain on a
finite state space can be improved in terms of asymptotic variance (or
in degenerate cases, not made worse) by transforming it to a Markov
chain on a related space that will be non-reversible (except when the
state space has only one or two states).

The non-reversible chains produced by this construction avoid, when
possible, transitions that backtrack by returning to the state from
which the chain just came.  This is done by expanding the state space
to pairs of states of the original chain --- representing, roughly
speaking, the previous and current states --- and then updating this
pair using two operations in sequence, one a swap, and the other a
modified Gibbs sampling update of the second component that tries to
avoid leaving the state unchanged.  Many such modifications are
possible; one that is generally applicable was introduced by Liu
(1996).  Though both the swap and the modified Gibbs sampling update
are reversible, their application in sequence is not reversible.

Simulation of this non-reversible chain will often require little or
no more computation time than simulation of the original reversible
chain.  The advantage of the non-reversible chain can be dramatic when
the original chain is such that suppressing backtracking has the
effect of forcing movement in the same direction for many steps,
thereby suppressing the slow random-walk motion that reversible chains
are subject to.  In other cases, however, the improvement may be
slight.  Aside from possible practical applications of the particular
construction I present, the results indicate that non-reversible
chains are fundamentally superior to reversible chains for MCMC
estimation, and hence research into improved MCMC methods may be best
directed toward methods based on non-reversible chains.

My proof that asymptotic variance will not increase as a result of
modifying the chain to avoid backtracking uses a new technique based on
dividing the chains into blocks delimited by transitions that are
affected by the modification, and then showing that the only effect of
the modification is to partially stratify the sampling for these
blocks.  This stratification can only decrease asymptotic variance, or
leave it unchanged.  As an introduction to this technique, I start by
giving a new proof of Peskun's (1973) theorem that asymptotic variance
will not be increased by modifying a reversible chain to reduce the
probability of staying in the same state, while keeping the
probability of other transitions at least as large as before.  This
proof gives some insight into why the hypothesis of reversibility is
necessary for Peskun's theorem.  This new proof technique holds
promise for proving that other transformations of both reversible and
non-reversible chains are also beneficial.

\section{\hspace*{-7pt}Preliminaries}\label{sec-prelim}\vspace*{-10pt}

Suppose we wish to estimate the expectation of some function, $f(x)$,
with respect to a distribution with probabilities $\pi(x)$, where $x$
is in some finite space $\X$.  (Generalizations to infinite spaces 
will not be dealt with in this paper.)  The MCMC approach to
this problem is to simulate a Markov chain $X_1,X_2,\ldots$ that has
$\pi$ as an invariant distribution --- that is, for which
\beq
   \pi(y) & = & \sum_{x\in\X} \pi(x) T(x,y),\ \ \ \mbox{for all $y\in\X$}
\label{eq-invariant}
\eeq
where $T(x,y) = P(X_{t+1}=y\,|\,X_t=x)$ are the transition probabilities
of the Markov chain (assumed to be the same for all $t$).  If the Markov 
chain is also irreducible (a series of transitions with non-zero probability
connects any two states), it will have only one invariant distribution, 
and the estimator\vspace*{-6pt}
\beq
   \hat\mu_n = {1 \over n} \sum_{t=1}^n f(X_t)
\label{eq-mcmc-est}
\eeq
will converge to $\mu = E_{\pi}[f(X)]$ as $n$ goes to infinity.
Furthermore, a Central Limit Theorem applies, showing that the distribution of 
$\hat\mu_n$ is asymptotically normal (possibly a
degenerate normal distribution with variance zero).
These fundamentals aspects of MCMC are discussed, for example, 
by Tierney (1994) and Liu (2001).  Some statements of the results mentioned
above in these references make a further assumption that the chain is 
aperiodic, but this is is not essential (Hoel, Port, and Stone 1972, Chapter 2,
Theorems 3, 5, and 7; Romanovsky 1970, Section 43).

The asymptotic variance of the estimator~(\ref{eq-mcmc-est}) is defined to be
\beq
   V_{\infty}(\hat\mu) & = & \lim_{n\rightarrow\infty} n\, \Var(\hat\mu_n)
\label{eq-asymvar}
\eeq
Note that this does not depend on the initial distribution for $X_1$.
The bias of the estimator will be of order $1/n$, so its asymptotic mean 
squared error will be equal to its asymptotic variance. In practice, rather 
than $\hat \mu_n$ from~(\ref{eq-mcmc-est}), we would use an estimator based 
only on $X_t$ with $t$ greater than some time past which we believe the chain
has reached a distribution close to $\pi$, but this refinement (which
reduces bias) does not affect the asymptotic variance.

Asymptotic variance can be used as a criterion for which of two Markov
chains with the same invariant distribution is better, on the
assumption that the squared error of a practical estimator based on $k$
consecutive states of the Markov chain will be approximately $k\,
V_{\infty}$.  This is not guaranteed to be true.  For example, if
$\pi$ is uniform over $\X$, a Markov chain that deterministically
cycles through all states in some order will have asymptotic variance
of zero, but if the number of states in $\X$ is enormous, an estimator
based on simulating a practical number of transitions of this Markov
chain may have large squared error.  Nevertheless, in many contexts,
$V_{\infty}$ will be a good guide to practical utility, and it will be
used in this paper as the criterion for comparing Markov chains.
Asymptotic variance has previously been used to compare Markov chains
by Peskun (1973) and by Mira and Geyer (2000), as well as by many
others.

A Markov chain is said to be ``reversible'' if its transition probabilities
satisfy the following ``detailed balance'' condition with respect to $\pi$:
\beq
      \pi(x)T(x,y) & = & \pi(y)T(y,x),\ \ \ \mbox{for all $x,y \in \X$}
\label{eq-dbal}
\eeq
As a consequence, a sequence $X_1,\ldots,X_n$ from a reversible Markov 
chain with $X_1$ having initial distribution $\pi$ will have the same 
distribution as the reversed sequence of states, $X_n,\ldots,X_1$.   
Detailed balance implies that $\pi$ is an invariant distribution of the
Markov chain, but the converse need not hold.

Many MCMC methods use reversible Markov chains, notably the
widely-used Metropolis-Hastings algorithm (Hastings 1970).  In this
algorithm, a transition from state $x$ is performed by first randomly
drawing a state, $x^*$, from some ``proposal distribution'', with
probabilities $S(x,x^*)$, and then accepting $x^*$ as the next state
of the chain with probability
\beq
   a(x,x^*) & = &
     \min \left[\, 1,\ { \pi(x^*)\, S(x^*,x) \over \pi(x)\ S(x,x^*) }\, \right]
\label{eq-MHaccept}
\eeq
If $x^*$ is not accepted, the next state of the chain is the same as the
current state.  The result is that for $y \ne x$, the transition probability
is $T(x,y) = S(x,y)a(x,y)$.
One can easily show that these transitions satisfy detailed balance with
respect to $\pi$, and hence leave $\pi$ invariant.  As a special case of
the Metropolis-Hastings algorithm, when $x$ consists of several components,
$x^*$ might differ from $x$ in only a single component, with the value
for that component in $x^*$ being drawn from its conditional distribution
under $\pi$ given the values of the other components.  The acceptance 
probability of~(\ref{eq-MHaccept}) will then alway be one.  This is called
a ``Gibbs sampling'' (or ``heatbath'') update of the component. 

Not all Markov chains used for MCMC are reversible, however.  In
particular, non-reversible Markov chains often arise as a result of
applying two or more reversible transitions in sequence.  If $T_1$ and
$T_2$ are matrices of transition probabilities that satisfy detailed
balance with respect to $\pi$ (and hence leave $\pi$ invariant), their
product, $T_1T_2$, will also leave $\pi$ invariant, but will typically
not satisfy detailed balance.  A common example of this is when Gibbs
sampling updates are applied to each component of state in some
deterministic sequence.

There is no reason to avoid non-reversible chains in practical
applications of MCMC --- what is essential is that the chain leave
$\pi$ invariant, not that it be reversible with respect to $\pi$.  The
non-reversibility of deterministic-scan Gibbs sampling is thought to
be of little significance, but other non-reversible MCMC methods are
designed to exploit non-reversibility to avoid the slow, diffusive
movement via a random walk that is typical of reversible Markov
chains.  Examples include ``overrelaxation'' methods (Adler 1981; Neal
1998, 2003) and the ``guided Monte Carlo'' methods of Horowitz (1991)
and Gustafson (1998). 

However, non-reversible Markov chains have often been avoided in
theoretical discussions, since they are harder to analyse than
reversible chains.  In contrast, Diaconis, Holmes, and Neal (2000)
analysed a particular non-reversible chain and showed that it
converges to its invariant distribution much faster than a related
reversible chain.  Mira and Geyer (2000) explored whether
non-reversible chains can be transformed to reversible chains with the
same asymptotic variance, and found a method that sometimes does this,
but not always, again showing that non-reversible chains might be
superior to reversible chains.  These results, and the practical
usefulness of some non-reversible chains, lead one to ask whether any
reversible chain can be transformed to a non-reversible chain that is
better.  With some caveats --- notably, a restriction to finite state
spaces --- this paper provides an affirmative answer to this question.

\section{\hspace*{-7pt}Peskun's theorem on modifying a reversible chain 
     to avoid \\
     \hspace*{-7pt}staying in the same state}\label{sec-peskun}\vspace*{-10pt}

Before showing how to construct a non-reversible chain that is better
than a given reversible chain, I will present Peskun's theorem, which
shows that modifying a reversible chain to decrease the probability of
staying in the same state, while keeping the probabilities of other
transitions at least as large, cannot increase asymptotic variance.
This theorem is relevant to the non-reversible construction that
follows.  I also introduce the new proof techniques I use by proving
this theorem in Section~\ref{sec-peskun-proof}.

{\em

\noindent \textbf{Theorem 1 (Peskun 1973):}
Let $X_1,X_2,\ldots$ and $X'_1,X'_2,\ldots$ be two irreducible
Markov chains on the finite state space $\X$, both of
which are reversible with respect to the distribution $\pi$, and hence have
$\pi$ as their unique invariant distribution.
Let the transition probabilities for these chains be 
\beq
T(x,y)=P(X_{t+1}=y\,|\,X_t=x),\ \ \ \ T'(x,y)=P(X'_{t+1}=y\,|\,X'_t=x)
\eeq
Let $f(x)$ be some function of state, whose expectation with respect to
$\pi$ is $\mu$.  Consider the following two estimators for $\mu$ based
on these two chains:
\beq
   \hat\mu_n = {1 \over n} \sum_{t=1}^n f(X_t),\ \ \ \ \
   \hat\mu'_n = {1 \over n} \sum_{t=1}^n f(X'_t)
\eeq
If $T$ and $T'$ satisfy the following condition,
\beq
   T'(x,y) \ \ge\ T(x,y),\ \ \ \mbox{for all $x,y \in \X$ with $x \ne y$}
\eeq
then the asymptotic variance of $\hat\mu'$ will be no greater than that
of $\hat\mu$.

}

Since it seems inefficient to stay in the same place, Peskun's theorem
might seem obvious.  Two facts show that the situation is more subtle than
this.
First, only the \textit{asymptotic} variance is guaranteed not to
increase if off-diagonal entries in the transition matrix are increased.
The variance of an estimator based on finite number of transitions,
started from $\pi$, may increase (Tierney, 1998).
Second, Peskun's theorem does not hold if the condition that the
chains be reversible is omitted.  Here is a counterexample using
a non-reversible chain with four states:\vspace*{5pt}

\hspace*{0.4in}\psfig{file=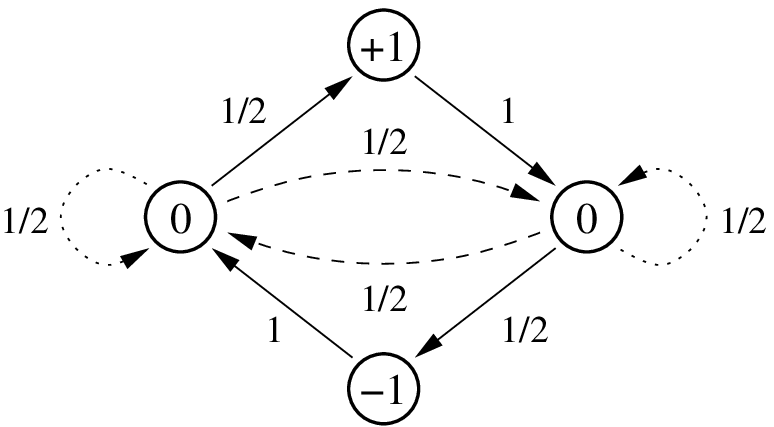,height=1.7in}

\vspace*{-1.65in}
\hfill\begin{minipage}[t]{2.4in}
 ~ \\[5pt]
 $\pi(x)\ =\,\ \scriptstyle  \left\{\!\!\begin{array}{ll} 
     \scriptstyle 1/3 & \!\mbox{\small where $ f(x)=0$} \\
     \scriptstyle 1/6 & \!\mbox{\small where $ f(x)\ne0$} \\[4pt]
                  \end{array}\right.$ \\
 ~ \\ 
 $\mu \ =\ 0$ 
\end{minipage}\\[52pt]
Values of $f(x)$ are shown in the circles.  Solid arrows show values
of both $T(x,y)$ and $T'(x,y)$; dotted arrows are for $T(x,y)$ only;
dashed arrows are for $T'(x,y)$ only.  The asymptotic variance is 
zero when using $T$ --- the chain proceeds clockwise through the
four states, never backtracking (though sometimes staying put in a 
state where $f(x)=\mu$), with the result that $|\hat\mu_n-\mu| \le 1/n$.
The modification that produces $T'$ disturbs this cyclic behaviour,
with the result that the asymptotic variance becomes greater than zero.


One application of Peskun's theorem is to motivate a modified form of
Gibbs sampling due to Liu (1996), which tries to avoid setting a component 
to the same value it had previously.  Suppose, for example, that the state 
consists of two components, and the current state is $(x,y)$.  As
mentioned above, a Gibbs sampling update for $y$ can be seen as a
Metropolis-Hastings update with a proposal distribution that keeps
the first component unchanged and draws a new value for the second component,
$y^*$, from its conditional distribution, $\pi(y|x)$.  In Liu's
modification, the proposal distribution is confined to values for the
second component other than the current value, $y$, with the probability 
for proposing $y^*$ being \mbox{$\pi(y^*|x)\,/\,(1\!-\!\pi(y|x))$}.  
The acceptance probability (from (\ref{eq-MHaccept})) then becomes
\beq
    a((x,y),\,(x,y^*))\ \ =\ \
    \min \left[\,1,\ 
     { \pi(x,y^*)\, \pi(y|x)\,/\,(1\!-\!\pi(y^*|x)) \over
       \pi(x,y)\, \pi(y^*|x)\,/\,(1\!-\!\pi(y|x))} \,\right]
    & = & 
    \min \left[\,1,\ { 1-\pi(y|x) \over 1-\pi(y^*|x) } \,\right]\ \ \ \
\label{eq-liuaccept}
\eeq
In the special case that $\pi(y|x)=1$, we never accept the proposal
(which is undefined).

One can easily verify that this modification increases
the probability of a transition to all values except
for the current value.  However, Peskun's theorem will not apply if
these modified Gibbs sampling updates are applied in sequence
(producing a non-reversible chain).  Peskun's theorem does apply if we
select a component to update at random, although this is not how Gibbs
sampling is commonly done in practice.  Liu's modification of Gibbs
sampling will play a role in the construction of a non-reversible
chain that avoids backtracking, which is presented next.

\section{\hspace*{-7pt}Constructing a non-reversible chain from a reversible 
     chain \\
     \hspace*{-7pt}so as to avoid backtracking}\label{sec-const}\vspace*{-10pt}

As above, suppose we have an irreducible Markov chain on a finite state space,
$\X$, with transition probabilities given by $T(x,y)$.  Suppose also
that this chain is reversible, so these transition probabilities
satisfy the detailed balance condition (\ref{eq-dbal}) with respect
to some invariant distribution, $\pi$.  In this section, I show how
to construct from $T$ a non-reversible Markov chain that avoids backtracking.
The state space for this chain will be $\ddot\X = \{ (x,y)\,:\, T(x,y)>0\}$, 
and it will leave invariant the distribution $\ddot\pi$ with probabilities 
defined as follows:
\beq 
    \ddot\pi(x,y) & = & \pi(x)T(x,y)\ \ =\ \ \pi(y)T(y,x)
\label{eq-ddotpi}
\eeq
One can view $\ddot\pi$ as the distribution for a pair of consecutive
states from the original chain, with the first state in the pair
drawn from $\pi$.  The second formula above follows from the reversibility
of the original chain.  Note that under $\ddot\pi$, the marginal distributions
of the first and second components are both $\pi$.  We can
therefore estimate the expectation with respect to $\pi$ of any function 
defined on $\X$ by averaging the values of either component in pairs 
distributed according to~$\ddot\pi$.

I will first show how to construct a chain with state space $\ddot\X$
and invariant distribution $\ddot\pi$ that is essentially the original
chain in disguise, but which can later be modified to prevent
backtracking.  This construction is called ``expanding'' the chain by
Kemeny and Snell (1960, Section 6.5).  A transition of this chain
consists of the following two operations, applied in sequence:\vspace*{-8pt}
\begin{enumerate}
\item[1)] Swap the two components of the state.\vspace*{-4pt}
\item[2)] Replace the second component of this swapped state with a 
          new value sampled from its conditional distribution 
          (under $\ddot\pi$) given the current value of the first 
          component.\vspace*{-8pt}
\end{enumerate}
From (\ref{eq-ddotpi}), the conditional probability for
the second component to be $y$, given that the first component is $x$,
is $T(x,y)$.  The transition probabilities, $\ddot T$, for this chain 
can therefore be written as follows:\vspace*{-9pt}
\beq
   \ddot T((x_0,y_0),\,(x_1,y_1)) & = & \delta(x_1,y_0)\,T(x_1,y_1)
\label{eq-ddotT}
\eeq
where $\delta(x,y)$ is one if $x=y$ and zero otherwise.

The first operation above leaves $\ddot\pi$ invariant, since
$\ddot\pi(x,y) = \ddot\pi(y,x)$, due to the reversibility of
the original chain.  The second operation above leaves $\ddot\pi$
invariant as well, since it is simply a Gibbs sampling update of the second
component.  Applying these two operations in sequence therefore also
leaves $\ddot\pi$ invariant.

Although both operations above are reversible, applying them in
sequence produces a non-reversible chain (except in degenerate
situations).  This non-reversibility is of no consequence, however,
since the above chain on $\ddot\X$ essentially replicates the
operation of the original chain on $\X$.  Starting from state
$(x_0,x_1)$, the chain will proceed to states $(x_1,x_2)$,
$(x_2,x_3)$, $(x_3,x_4)$ etc., with each $x_t$ being drawn according
to the probabilities $T(x_{t-1},x_t)$, just as in the original chain.
If we estimate the expectation of $f$ with respect to $\pi$ by the
average value of $f$ applied to the second components of these states,
the result will be exactly the same as an estimate based on states of 
the original chain.

To obtain a more interesting non-reversible chain, we can change the
second operation above to use the modified Gibbs sampling update of
Liu (1996), discussed above in Section~\ref{sec-peskun}.  More
generally, we might modify the second operation in any way that
reduces the probability of staying in the same state, while keeping
the probabilities of transitions to other states at least as large as
before, and maintaining reversibility with respect to $\ddot\pi$.
Such a modified chain, whose transition probabilities we will write as
$\ddot T'$, will differ substantively from the original reversible
chain, since it will reduce the probability of ``backtracking'' to the
state preceding the current state.  For example, starting from state
$(x_0,x_1)$, a chain with transition probabilities $\ddot T$,
equivalent to the original reversible chain, might proceed to state
$(x_1,x_2)$ and then to $(x_2,x_1)$ --- corresponding to the original
reversible chain moving from $x_1$ to $x_2$ and then back to $x_1$.  A
modified chain with transitions $\ddot T'$ might also proceed from
$(x_0,x_1)$ to $(x_1,x_2)$, but after the swap operation of the next
transition, the state $(x_2,x_1)$ would be updated by a modified Gibbs
sampling operation that has a reduced probability of leaving the
second component equal to $x_1$.

That this avoidance of backtracking cannot increase asymptotic
variance is the central result of this paper.  This is stated in
the theorem below, which is proved in Section~\ref{sec-main-proof}.

{\em

\noindent \textbf{Theorem 2:}
Let $X_1,X_2,\ldots$ be an irreducible
Markov chain on the finite state space $\X$ having
transition probabilities $T(x,y)=P(X_{t+1}=y\,|\,X_t=x)$ that
satisfy detailed balance with respect to the distribution with
probabilities $\pi(x)$.
Define a Markov chain $(X_1,Y_1),\,(X_2,Y_2),\,\ldots$ on the
state space $\ddot\X = \{ (x,y)\,:\, T(x,y)>0\}$ with transition probabilities
\beq
   \ddot T'((x_0,x_1),\,(y_0,y_1)) & = & \delta(x_1,y_0)\,U'_{x_1}(x_0,y_1)
\label{eq-ddotTprime}
\eeq
where $\delta(x,y)$ is one if $x=y$ and zero otherwise, and $U'_x(y,z)$
defines a set of probabilities for $z\in\X$ for any values of $x,y\in\X$,
satisfying the following two conditions for all $x,y,z \in \X$ with $y \ne z$:
\beq
   T(x,y)\,U'_x(y,z) & = & T(x,z)\,U'_x(z,y) 
\label{eq-Ucond1}\\[4pt]
   U'_x(y,z) & \ge & T(x,z)
\label{eq-Ucond2}
\eeq
Let $f(x)$ be some function of state, whose expectation with respect to
$\pi$ is $\mu$.  Define the following two estimators for $\mu$ based
on these two chains:\vspace*{-3pt}
\beq
   \hat\mu_n = {1 \over n} \sum_{t=1}^n f(X_t),\ \ \ \ \
   \hat\mu'_n = {1 \over n} \sum_{t=1}^n f(Y'_t)
\eeq
Then the following properties of the 
$\ddot T'$ and the estimators above hold:\ \ (a) the chain with transition 
probabilities $\ddot T'$ is irreducible; (b) the transition probabilities 
$\ddot T'$ leave invariant the distribution with probabilities 
$\ddot\pi(x,y)=\pi(x)T(x,y)$; (c) if $\X$ contains at least three elements, 
the chain with transition probabilities $\ddot T'$ is not reversible with 
respect to $\ddot\pi$; (d) the bias of the estimator 
$\hat\mu'_n$ is of order $1/n$; (e) the asymptotic variance of $\hat\mu'$ 
is no greater than the asymptotic variance of $\hat\mu$. \\[6pt]
}
The transition probabilities $\ddot T$ (from~(\ref{eq-ddotT})), which 
essentially mimic $T$, can also be written as $\ddot T((x_0,x_1),(y_0,y_1)) = 
\delta(x_1,y_0)\,U_{x_1}(x_0,y_1)$, with $U_{x_1}(x_0,y_1)=T(x_1,y_1)$.  
The change from $T$ to $\ddot T'$
in this theorem can therefore also be seen as a change from $\ddot T$
to $\ddot T'$ or from $U$ to $U'$.

The new probabilities $U'_x(y,z)$ are modified update probabilities for 
the second component of state,
with $x$ being the first component of state, $y$ the
current value of the second component, and $z$ a new value for the
second component.  These updates must satisfy detailed
balance with respect to $\ddot\pi$.  In the case of
Liu's modification, we find using (\ref{eq-liuaccept}) that
for all $x,y,z\in\X$ with $z \ne y$,
\beq
   U'_x(y,z) & = & { T(x,z) \over 1-T(x,y) } \
                  \min\left[\, 1,\ { 1-T(x,y) \over 1-T(x,z)} \,\right]
          \ \ =\ \ \min\left[\, { T(x,z) \over 1-T(x,y) },\ 
                                { T(x,z) \over 1-T(x,z)} \,\right]\ \
\label{eq-liumod}
\eeq
$U'_x(y,y)$ is determined from the above by the requirement
that probabilities sum to one.  Note that if $T(x,y) \ge 1/2$, this
expression simplifies to $T(x,z)\,/\,(1\!-\!T(x,z))$.

As a first example of such a modified chain, consider a Markov chain
on $\X=\{ 1,\,2,\,\ldots,\,N \}$ with transition probabilities of
$T(x,y)=1/2$ when $y=x+1$ or $y=x-1$ or $x=y=0$ or $x=y=N$, and
$T(x,y)= 0$ otherwise.  This chain is irreducible, has the uniform
distribution as its invariant distribution, and is reversible.
From~(\ref{eq-liumod}), we can see that for given $(x,y)\in\ddot\X$,
$U'_x(y,z)=1$ for some $z$.  The transitions, $\ddot T'$, of 
the modified chain are therefore deterministic.  For $N=5$, these
transitions follow the arrows in the the diagram below:\\[12pt]
  \centerline{\psfig{file=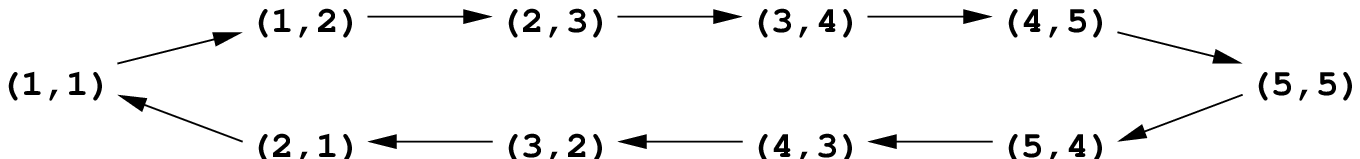}}\\[12pt]
The periodic nature of the modified chain results in the asymptotic variance 
being zero for any function of state, whereas for the original chain, the 
asymptotic variance is of order $N^2$, due to its random walk behaviour.
This example parallels the chain analysed by Diaconis, Holmes, and Neal
(2000), with $c$ set to zero in the definition (their equation 4.1) of 
their chain.  Note, however, that the more general scheme they 
describe (in their section~5.1) does not correspond to the result of 
modifying a random walk Metropolis algorithm to avoid backtracking in 
the way described here.

As another illustration, consider a chain on $\X = \{
1,\,2,\,\ldots,\,N\} \times \{ 1,\,2,\,\ldots,\,M\}$, which may be
visualized as dots arranged in an $N$ by $M$ rectangle, in which
transitions go up, down, left, or right, with equal probabilities,
except that if such a movement would leave the rectangle, the chain
instead stays in the current state.  This chain leaves the uniform 
distribution invariant.  These transitions
are shown below, for $N=6$ and $M=3$ (the unmarked transition probabilities
are $1/4$):\\[12pt]
  \centerline{\psfig{file=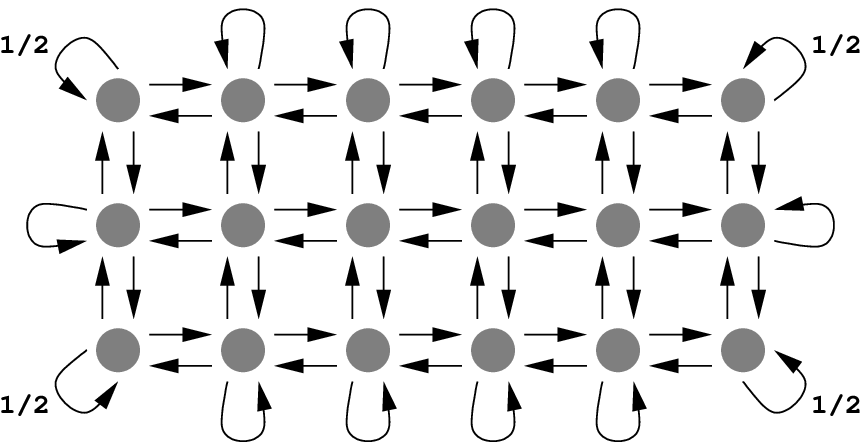}}

The state space, $\ddot\X$, of the modified chain consists of the arrows
in the diagram above.  The transitions probabilities, $\ddot T'$, for
the modified chain, based on the modified updates of~(\ref{eq-liumod})
are illustrated below, for two possible current states:\\[12pt]
  \centerline{\psfig{file=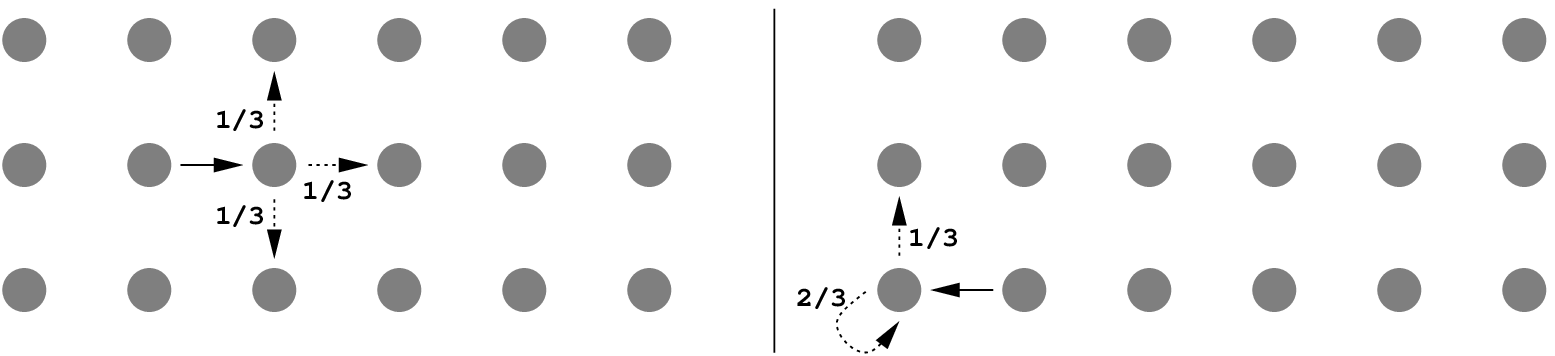}}\\[6pt]
The current states in these diagrams are shown by sold arrows, and
possible successor states by dotted arrows, labeled with their
probabilities.  

The diagrams below show two paths within the rectangle,
produced using the original chain (on the left) and the modified chain
(on the right):\\[8pt]
  \centerline{\psfig{file=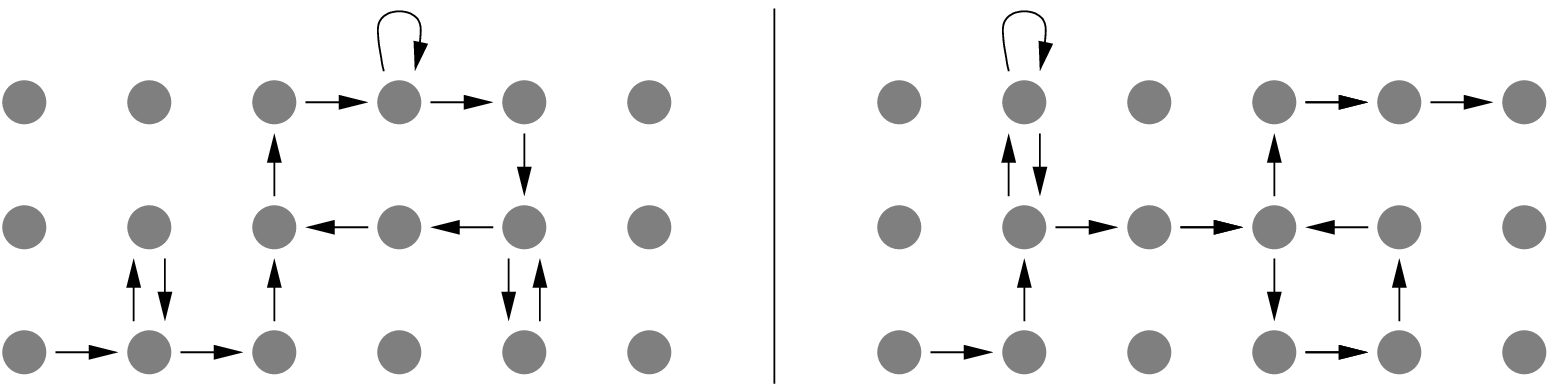}}\\[12pt]
Note that in two places the original chain backtracks to the preceding
state.  The new chain never backtracks in this way, but it is still
possible for it to revisit states that were visited two or more time
steps earlier.  As a result, the improvement in asymptotic variance
is not as dramatic as for the previous example.  Asymptotic variance
is improved only by a constant factor, which does not increase with $N$ and 
$M$.

To simulate a chain that has been modified to avoid backtracking, with
transition probabilities $\ddot T'$, we need to be able to draw a
value from $\X$ according to the probabilities $U'_x(y,\,\cdot\,)$.  If
we use $U'_x(y,z)$ defined by (\ref{eq-liumod}), and if $T(x,z)$ is
non-zero for only a small, known set of $z$ values, we can do this by
explicitly computing the probabilities using (\ref{eq-liumod}).  This
will often be about as efficient as simulating the original chain.

When $T(x,z)$ is non-zero for many value of $z$, the following
procedure for drawing a $z$ value from $U'_x(y,\,\cdot\,)$ as defined
by~(\ref{eq-liumod}) may be useful.  If $T(x,y) \ge 1/2$, draw a value
$z^*$ according to the probabilities $T(x,z^*)$.  If $z^*=y$, let
$z=y$.  Otherwise, accept $z^*$ as the value $z$ with probability
$1/(1\!-\!T(x,z^*))$.  If $z^*$ is not accepted, let $z$ equal $y$.
If instead $T(x,y) < 1/2$, repeatedly draw $z^*$ according to the
probabilities $T(x,z^*)$, until a $z^*$ not equal to $y$ is obtained
(which won't take long).  Accept this $z^*$ as the value $z$ with
probability $\min[\, 1,\ (1\!-\!T(x,y))\,/\,(1\!-\!T(x,z)\,]$.  If $z^*$
is not accepted, let $z$ equal $y$.

In some cases, neither of the two procedures described above for
simulating from $U'_x(y,\,\cdot\,)$ may be easy to implement
efficiently --- for instance, $T(x,y)$ for a Metropolis-Hastings
transition may be hard to compute when $y=x$, since this requires
summing the probabilities of rejection for all possible proposals.  In
any case, a rigorous demonstration that a modified chain that avoids
backtracking can be simulated as quickly as the original chain is too
much to expect, because it might be possible to simulate the original chain
especially easily using some special trick that is not applicable to the 
modified chain.

Nevertheless, I think it is fair to say that avoiding backtracking,
either using Liu's modified Gibbs sampling update or some other form
for $U'_x(y,z)$, is not the sort of modification that inherently
involve a large increase in computation time per transition.  That
this modification decreases asymptotic variance (or in degenerate
cases, does not increase it) is therefore an important indication that
non-reversible chains have an advantage over reversible chains.

\section{\hspace*{-7pt}A new proof of Peskun's 
         theorem}\label{sec-peskun-proof}\vspace*{-10pt}

As an introduction to the techniques that will be used to prove that
the no-backtracking construction of the previous section does not
increase asymptotic variance, I will here use these techniques to
prove Peskun's theorem, stated as Theorem~1 in Section~\ref{sec-peskun}.

In this proof, the ``old chain'' will refer to the original chain with
transition probabilities $T$, and the ``new chain'' will refer to the
chain with transition probabilities $T'$, which may be smaller than
those of the old chain for self transitions, but are at least as large
for transitions between distinct states.  The proof that the estimator
for the expectation of any function of state using the new chain has
asymptotic variance at least as small as the corresponding estimator
using the old chain will proceed as follows:\vspace*{-8pt}\vfil
\begin{enumerate}

\item[1)] We reduce the problem to comparing asymptotic 
      variances when $T$ and $T'$ differ only for transitions involving 
      two states, $A$ and $B$.\vspace*{-4pt}\vfil

\item[2)] We can view simulations of the old and new chains as differing 
      only for certain ``delta'' transitions involving states $A$ and 
      $B$.\vspace*{-4pt}\vfil

\item[3)] These delta transitions divide the Markov chain simulation into blocks
      of states, which start and end in either state $A$ or state $B$. 
      We can rewrite the old and new estimators, $\hat\mu$ and $\hat\mu'$, 
      in terms of the lengths of these blocks and the sums of the function
      values for states in these blocks.\vspace*{-4pt}\vfil

\item[4)] We see that blocks starting and ending with $A$ and blocks 
      starting and ending with $B$ are equally likely, but may have
      different distributions for their contents.
      In contrast, blocks that start with $A$ and end with $B$ have essentially
      the same distribution of content as blocks that start with $B$ and end
      with $A$.\vspace*{-4pt}\vfil

\item[5)] The only difference between the old and new chains is that in the
      new chain the sampling for ``homogeneous'' blocks (starting and ending
      in the same state) is \textit{stratified} --- 
      there are the \textit{same} number of blocks starting and ending with 
      $A$ as blocks starting and ending with $B$, whereas the split
      between these types is random in the old chain (albeit with equal
      probabilities for the two types of homogeneous blocks).\vspace*{-4pt}\vfil
      
\item[6)] Finally, this stratification will lower (or at least not 
      increase) the asymptotic variance.\vspace*{-8pt}\vfil

\end{enumerate}

\subsection*{Step 1:\ \ Looking at one pair of states is enough}\vspace*{-7pt}

Whenever $T'(x,y) \ge T(x,y)$ for all $x \ne y$, we can get from $T$ to
$T'$ by a series of steps that each change transition probabilities for only
a single pair of states.  For example, consider the following steps from
$T$ to $T'$ that both satisfy detailed balance with respect to 
$\pi = [\, 0.4\ \ 0.4\ \ 0.2 \,]$:
\beq
  T\ =\ 
  \left[\begin{array}{lll} 0.4 & 0.4 & 0.2 \\ 0.4 & 0.4 & 0.2 \\ 0.4 & 0.4 & 0.2
              \end{array}\right]
  \Rightarrow 
  \left[\begin{array}{lll} 0.3 & 0.5 & 0.2 \\ 0.5 & 0.3 & 0.2 \\ 0.4 & 0.4 & 0.2
              \end{array}\right]
  \Rightarrow 
  \left[\begin{array}{lll} 0.3 & 0.5 & 0.2 \\ 0.5 & 0.2 & 0.3 \\ 0.4 & 0.6 & 0.0
  \end{array}\right]\ =\ T'
\eeq
Furthermore, if the detailed balance condition~(\ref{eq-dbal}) holds for
$T$ and for $T'$, it will hold also for all the intermediate transition 
probabilities (such as those in the middle matrix above), since the pair 
of transition probabilities for any
$x$ and $y$ (with $x \ne y$) at any intermediate point will be either
the same as for $T$ or the same as for $T'$.

It is therefore enough to prove Peskun's Theorem when $T$ and $T'$
differ for only two states, say $A$ and $B$.  The transition probabilities
for the old and new chain will then be related as follows:
\beq
  T'(x,y)\ =\ T(x,y),\ \ \ 
     \mbox{when $x \notin \{A,B\}$ or $y \notin \{A,B\}$~~~~} 
\nonumber\\
\label{eq-Trels}\\[-10pt]
    \begin{array}{ll} 
        T'(A,A)\ =\ T(A,A)-\delta_A,\ \ & T'(A,B)\ =\ T(A,B)+\delta_A \\
        T'(B,A)\ =\ T(B,A)+\delta_B,\ \ & T'(B,B)\ =\ T(B,B)-\delta_B
    \end{array}
\nonumber
\eeq
where $\delta_A$ and $\delta_B$ are positive.

\subsection*{Step 2:\ \ Marking ``delta'' transitions}\vspace*{-7pt}

Transitions $T$ and $T'$ differ only if the current state is $A$ or
$B$, and then only with respect to how a probability mass of
$\delta_A$ or $\delta_B$ is assigned to new states $A$ or $B$.  We can
mark such ``delta'' transitions while simulating the Markov chain.

The standard way to simulate a Markov chain is as follows:\ \ For each state, 
$x$, partition the interval $[0,1)$ into intervals $[\ell(x,y),h(x,y))$ such 
that $h(x,y)-\ell(x,y)\ =\ T(x,y)$; to simulate a transition out of state
$x$, generate a random variate, $U$, that is uniformly distributed on 
$[0,1)$, and move to the state, $y$, for which $\ell(x,y) \le U < h(x,y)$.
We can choose to simulate the old transitions, $T$, using partitions in which
$\ell(A,A) = \ell(B,B) = 0$.  With such a choice, we can write the algorithm
for simulating a transition of the old chain in the manner on the left below, 
in which a slight change yields the simulation algorithm for the new chain 
shown below on the right:\vspace*{7pt}

{\small\sf

\hspace*{20pt}\begin{minipage}{3in}
\textbf{Old chain:} \\[3pt]
$U \sim \mbox{Uniform}(0,1)$ \\
if $X_t = A$ and $U<\delta_A$ then \\
\hspace*{9pt}$X_{t+1} = A$,\ \ mark this transition \\
else if $X_t = B$ and $U<\delta_B$ then \\
\hspace*{9pt}$X_{t+1} = B$,\ \ mark this transition \\
else \\
\hspace*{9pt}$X_{t+1} = y$\ such that $U \in [\ell(X_t,y),h(X_t,y))$
\end{minipage}\hspace*{9pt}%
\begin{minipage}{3in}
\textbf{New chain:} \\[3pt]
$U \sim \mbox{Uniform}(0,1)$ \\
if $X_t = A$ and $U<\delta_A$ then \\
\hspace*{9pt}$X_{t+1} = B$,\ \ mark this transition \\
else if $X_t = B$ and $U<\delta_B$ then \\
\hspace*{9pt}$X_{t+1} = A$,\ \ mark this transition \\
else \\
\hspace*{9pt}$X_{t+1} = y$\ such that $U \in [\ell(X_t,y),h(X_t,y))$
\end{minipage}

}

\vspace*{6pt}

\noindent Clearly, $T$ and $T'$ differ only for the ``delta'' transitions 
marked above.

\subsection*{Step 3:\ \ Using delta transitions to define blocks}\vspace*{-7pt}

We can use the markings of delta transitions to divide a simulation of
one of these Markov chains into ``blocks'' of consecutive states, that
both start and end with either state $A$ or state $B$.  Note that
states $A$ and $B$ may also occur at places other than the start and
end of a block.  It is possible for a blocks to consist of only a
single $A$ or a single $B$.

Since asymptotic variance does not depend on the initial state
distribution, let's suppose that\linebreak $P(X_1=A) \, =\, P(X_1=B)
\,=\, 1/2$, so that the chains will begin at the start of a block.

For the old chain, with transitions $T$, we might see blocks like 
this:\vspace*{5pt}

\hspace*{0.8in}\psfig{file=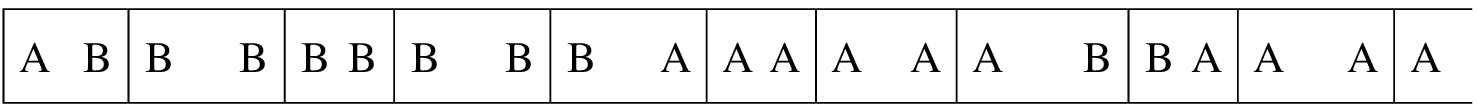,height=0.32in}

\noindent For the new chain, with transitions $T'$, the blocks might 
look like this:\vspace*{5pt}

\hspace*{0.8in}\psfig{file=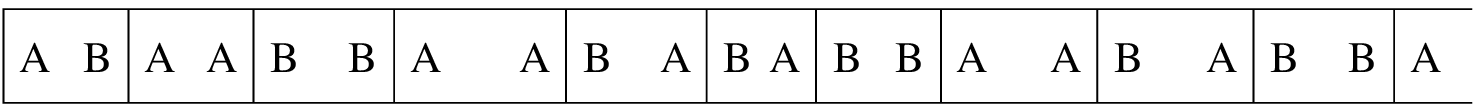,height=0.32in}

\noindent The difference is that in the old chain, the state stays the same
when crossing a block boundary, whereas for the new chain, it changes
from $A$ to $B$ or from $B$ to $A$.

We can view the simulation in terms of these blocks, and write the
estimates $\hat\mu$ and $\hat\mu'$ in terms of the lengths of the blocks
and the sums of $f$ for states in these blocks.  For the old chain,
\beq
  \hat\mu_n & \approx & \sum_{i=1}^k H_i \,\Big/\, \sum_{i=1}^k L_i
\eeq
where $H_i$ is the sum for $f(X_t)$ for states $X_t$ in block $i$,
$L_i$ is the length of block $i$, and $k$ is the number of 
blocks in the $n$ iterations of the chain.  The equality is only
approximate because there may be a partial block after block $k$.
Estimation in terms of blocks is discussed further in Step~6 of the proof.

\subsection*{Step 4: 
  Probabilities of the four types of blocks and their contents}\vspace*{-7pt}

Blocks come in four types --- $AA$, $BB$, $AB$, $BA$ --- based on 
start and end states.  For both the old and new chains, the 
probabilities of these types (ie, their frequencies of occurrence in a long
realization of the chain) satisfy
\beq
   P(AA) = P(BB)\ \ \mbox{and}\ \ P(AB) = P(BA) 
\eeq

We can show this using the fact that both $T$ and $T'$ leave $\pi$
invariant.  In particular, for the old chain,\vspace*{-10pt}
\beq
  \pi(B) & = & \pi(A)T(A,B)\ +\ \pi(B)T(B,B) \ +\!\!\! 
               \sum_{x \notin \{A,B\}}\!\!\! \pi(x)T(x,B) \\[-22pt]\nonumber
\eeq
while for the new chain, using the relationships in (\ref{eq-Trels}),
\beq
  \pi(B) & = & \pi(A)(T(A,B)+\delta_A)\ +\ \pi(B)(T(B,B)-\delta_B) \ +\!\!\! 
               \sum_{x \notin \{A,B\}}\!\!\! \pi(x)T(x,B) \\[-22pt]\nonumber
\eeq
from which it follows that $\pi(A)\delta_A\ =\ \pi(B)\delta_B$.

This lets us show that for a state, $X_t$, from the old chain (with $t$ 
being large),
\beq
   P(X_t\ \mbox{starts block with $A$}) 
    \,=\, 
   P(X_{t-1}\!=\!A)\,P(\mbox{delta transition at $t\!-\!1$}\,|\,X_{t-1}\!=\!A)
   & \!\!=\!\! & \pi(A)\,\delta_A \ \ \ \ \ \ \ \\
   P(X_t\ \mbox{starts block with $B$}) 
    \,=\, 
   P(X_{t-1}\!=\!B)\,P(\mbox{delta transition at $t\!-\!1$}\,|\,X_{t-1}\!=\!B)
   & \!\!=\!\! & \pi(B)\,\delta_B\ \ \ \ \ \ \
\eeq
and hence $P(X_t\ \mbox{starts block with $A$}) 
\ =\ P(X_t\ \mbox{starts block with $B$})$.  In the same way, we
see that $P(X_t\ \mbox{ends block with $A$}) 
\ =\ P(X_t\ \mbox{ends block with $B$})$.  It follows that
\beq
   P(AA) + P(AB)\ =\ P(BB) + P(BA)\ \ \mbox{and}\ \
   P(AA) + P(BA)\ =\ P(BB) + P(AB)
\eeq
so $P(AA) = P(BB)$ and $P(AB) = P(BA)$. 

Similarly, for a state, $X_t$, from the new chain (with $t$ 
being large),
\beq
   P(X_t\ \mbox{starts block with $A$}) 
    \,=\, 
   P(X_{t-1}\!=\!B)\,P(\mbox{delta transition at $t\!-\!1$}\,|\,X_{t-1}\!=\!B)
   & \!\!=\!\! & \pi(B)\,\delta_B \ \ \ \ \ \ \ \\
   P(X_t\ \mbox{starts block with $B$}) 
    \,=\, 
   P(X_{t-1}\!=\!A)\,P(\mbox{delta transition at $t\!-\!1$}\,|\,X_{t-1}\!=\!A)
   & \!\!=\!\! & \pi(A)\,\delta_A\ \ \ \ \ \ \
\eeq
and hence $P(X_t\ \mbox{starts block with $A$}) 
\ =\ P(X_t\ \mbox{starts block with $B$})$ for the new chain as well,
and similarly $P(X_t\ \mbox{ends block with $A$}) \ =\ P(X_t\ 
\mbox{ends block with $B$})$,
from which it again follows that $P(AA) = P(BB)$ and $P(AB) = P(BA)$. 

Although blocks of type AA and blocks of type BB are equally common,
the distributions for their contents --- and hence for their length
and for the sum of values of $f(x)$ over states in the block --- will
generally be different. In contrast, blocks of type AB and blocks of
type BA have the \textit{same} distribution of content --- except that
the BA blocks are the reversals of the AB blocks, which has no effect
on the sum of $f(x)$ for states in the block.  This equivalence of 
$AB$ and $BA$ blocks is a consequence of the chains being reversible, 
and holds for both the old and new chains.

To illustrate:\ \ The probability of block $AQB$ occurring at
some large time $t$ in the old chain is
\beq
\lefteqn{P(X_t=A\ \mbox{\& block starts})\,P(X_{t+1}=Q\,|\,X_t=A)
   P(X_{t+2}=B\ \mbox{\& block ends}\,|\,X_{t+1}=Q)}\ \ \ \nonumber\\
   & = &
   \pi(A)\,\delta_A\, T(A,Q)\, T(Q,B)\, \delta_B \ =\
   \delta_A\delta_B\, \pi(A)\, T(A,Q)\,T(Q,B)
   \ \ \ \ \ \ \ \ \ \ \ \ \ \ \ \ \ \ \ \ \ \ \ \ \ \ \ \ \\
   & = &
   \delta_A\delta_B\, T(Q,A)\, \pi(Q)\, T(Q,B) \ =\ 
   \delta_A\delta_B\, T(Q,A)\, T(B,Q)\, \pi(B) \\ 
   & = & 
   \pi(B)\,\delta_B\, T(B,Q)\, T(Q,A)\, \delta_A
\eeq
which is also the probability of block $BQA$ occurring at time $t$.
For the new chain, the probability of block $AQB$ occurring at time $t$ is
\beq
\lefteqn{P(X_t=A\ \mbox{\& block starts})\,P(X_{t+1}=Q\,|\,X_t=A)
   P(X_{t+2}=B\ \mbox{\& block ends}\,|\,X_{t+1}=Q)}\ \ \ \nonumber\\
   & = &
   \pi(B)\,\delta_B\, T(A,Q)\, T(Q,B)\, \delta_B \ =\
   \pi(A)\,\delta_A\, T(A,Q)\, T(Q,B)\, \delta_B
   \ \ \ \ \ \ \ \ \ \ \ \ \ \ \ \ \ \ \ \ \ \ \ \ \ \ \ \ 
\eeq
which is the same as for the old chain, and the same as for block $BQA$.

\subsection*{Step 5: In the new chain, sampling for homogeneous blocks 
             is stratified}\vspace*{-7pt}

Rather than simulate the chains one state at a time, let's imagine
simulating the chain block by block.  To show the relationship between
the old and the new chains, I'll show how this simulation can be done
in a coupled fashion.

To do this, we will need 
the probability that a block is ``homogeneous'' --- that it ends with the same
state it begins with --- which is\vspace*{-5pt}
\beq
  P(\mbox{ends with $A$}\,|\,\mbox{starts with $A$}) & \!=\! &
  {P(AA) \over P(AA)+P(AB)} \nonumber\\
  & \!=\! & {P(BB) \over P(BB)+P(BA)} 
  \ = \ P(\mbox{ends with $B$}\,|\,\mbox{starts with $B$})\ \ \ \ \ \
\eeq
Call this probability $h$, and note that it is the same for the old chain
and the new chain, since the transitions within a block, and the marking
of its end, are the same for both chains.  Note as well that the distribution
of the contents of a block, given its type, is the same for the old chain
and the new chain.

We can now simulate block transitions for the ``old'' and ``new'' chains as 
follows.  We'll assume $H$ below is sampled the same for both chains,
but that the simulation of the contents of blocks is not coupled
between the old and new chains.\vspace{4pt}

{\small\sf

\hspace*{20pt}\begin{minipage}{3in}
\textbf{Old chain:} \\[1pt]
$H \sim \mbox{Bernoulli}(h)$ \\
if $H=1$ then \\
\hspace*{10pt}if previous block ended with $A$ \\
\hspace*{20pt}simulate an $AA$ block \\%
\hspace*{10pt}else \\
\hspace*{20pt}simulate a $BB$ block \\
else \\
\hspace*{10pt}if previous block ended with $A$ \\
\hspace*{20pt}simulate an $AB$ block\\%
\hspace*{10pt}else \\
\hspace*{20pt}simulate an $AB$ block, then reverse it
\end{minipage}\hspace*{10pt}%
\begin{minipage}{3in}
\textbf{New chain:} \\[1pt]
$H \sim \mbox{Bernoulli}(h)$ \\
if $H=1$ then \\
\hspace*{10pt}if previous block ended with $A$ \\
\hspace*{20pt}simulate a $BB$ block \\%
\hspace*{10pt}else \\
\hspace*{20pt}simulate an $AA$ block \\
else \\
\hspace*{10pt}if previous block ended with $A$ \\
\hspace*{20pt}simulated an $AB$ block, then reverse it\\%
\hspace*{10pt}else \\
\hspace*{20pt}simulate an $AB$ block
\end{minipage}

}

\vspace*{6pt}

Comparing the simulations for the old and new chains, we see that they
produce the \textit{same} sequence of homogeneous/non-homogeneous
blocks.  However, for the new chain, the homogeneous blocks
\textit{alternate} between $AA$ blocks and $BB$ blocks.  
This is true both when one homogeneous block follows another, and when
any number of non-homogeneous blocks intervene.  In the old chain, the
type of homogeneous block changes only when an odd number of
non-homogeneous blocks intervene.  This is illustrated below:\vspace{3pt}

\hspace*{0.2in}Old chain:\vspace*{-19pt}

\hspace*{1.05in}\psfig{file=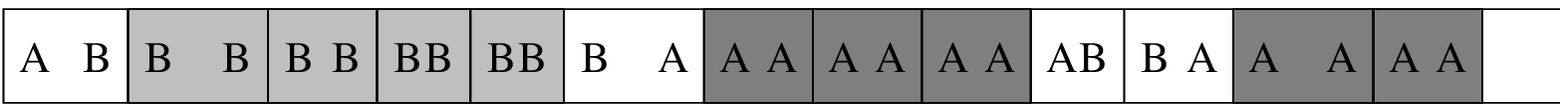,height=0.32in}

\hspace*{0.2in}New chain:\vspace*{-19pt}

\hspace*{1.05in}\psfig{file=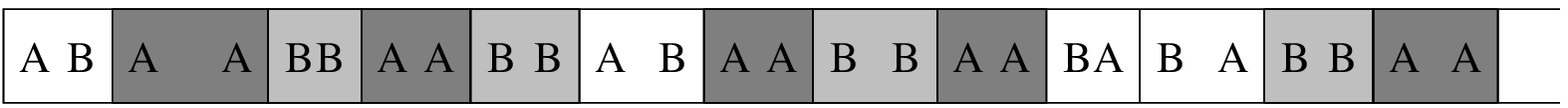,height=0.32in}

Because $AA$ blocks alternate with $BB$ blocks in the new chain,
sampling within the new chain is \textit{stratified} in this respect
--- that is, the number of $AA$ blocks will be equal to the number of
$BB$ blocks (plus or minus one).  We can also see this by noting that
in the new chain, every block ending in $A$ (except the last) is paired
with a following block beginning with $B$.  Letting $N_{AA}$, $N_{AB}$,
$N_{BA}$, and $N_{BB}$ be the numbers of blocks of each type, it 
follows that 
\beq
   |\, (N_{AA}\!+\!N_{BA})\, -\, (N_{BB}+N_{BA})\,| & \le & 1
\eeq
and hence $|N_{AA}-N_{BB}| \le 1$.

\subsection*{Step 6: Stratification of homogeneous blocks won't
             increase asymptotic variance}\vspace*{-7pt}

The intuition behind the proof is now complete: Sampling with the new
chain is stratified with respect to blocks of type $AA$ and $BB$.
Furthermore, this is the \textit{only} difference between the old and
new chains, since the sum of the values of $f$ for states in a block has 
the same distribution for $AB$ blocks and $BA$ blocks.  Since we expect that
stratification will not increase asymptotic variance (and will
typically decrease it), the asymptotic variance for the new chain
should be no larger than for the old chain.

To justify this formally, we need two lemmas whose detailed
statements, proofs, and applications to this proof are found in the
Appendix.  The first lemma says that the asymptotic variance
(appropriately defined) of an estimator based on a simulation that
continues for some specified number of blocks is the same as that of
an estimator based on a simulation that continues for some specified
number of Markov chain transitions.  This is true because the Central
Limit Theorem implies that asymptotically there is very little
difference between simulating for a specified number of blocks and
simulating for the number of transitions equal to the expected number
for that many blocks.

Accordingly, we can compare the old and new chains in the context of
simulations that continue for a specified number of blocks.  The
second lemma justifies the idea that the stratification of $AA$ and
$BB$ blocks will not increase the asymptotic variance of estimators
based on such simulations.  Stratification is of course well-known to
be beneficial (or at least not harmful) in the context of independent
sampling from two populations.  The lemma shows that this continues to
be true when, as here, the stratification is only partial (the ratio
of homogeneous to non-homogeneous blocks is not fixed), sampling has a
Markov chain aspect rather than being independent, and the estimator
takes the form of a ratio rather than a linear function of the sampled
variables.

\vspace*{10pt}

This proof provides some insight into why Peskun's theorem needs the
premise that the chains are reversible with respect to $\pi$, rather
than the weaker premise that they leave $\pi$ invariant.  This premise
is used at two points in the proof.  In Step~1, the reduction to old
and new chains that differ only for transitions involving two states
would not be possible for non-reversible chains, since there would be
no guarantee that the intermediate chains linking old and new chains
differing for several pairs of states would leave $\pi$ invariant.
This is not relevant to the counterexample in
Section~\ref{sec-peskun}, however, since it involves two chains that
already differ only with regard to transitions between two states 
(state $A$ on the left and $B$ on the right).

The reason the new chain in the counterexample has higher asymptotic
variance relates to the second use of reversibility, in Step 4, where
the contents of blocks of type $BA$ are seen to have essentially the
same distribution as the contents of blocks of type $AB$, apart from a
reversal that does not affect the sum of $f(x)$, which is what matters
for the estimates.  For the counterexample, the sum of $f(x)$ for
blocks of type $AB$ will always be $+1$, whereas this sum will always
be $-1$ for blocks of of type $BA$.  In contrast, the sum of $f(x)$
for blocks of type $AA$ or type $BB$ will always be $0$.  Examining
Step~5 of the proof, one can see that while the new chain stratifies
sampling for blocks of type $AA$ and $BB$, the old chain stratifies
sampling for blocks of type $AB$ and $BA$.  For a non-reversible
chain, stratifying between types $AB$ and $BA$ may be more important,
so it is possible for the old chain to have lower asymptotic variance
than the new chain.

\section{\hspace*{-7pt}Proof that modifying a reversible chain to avoid
    backtracking \\
    \hspace*{-7pt}doesn't increase asymptotic 
    variance}\label{sec-main-proof}\vspace*{-10pt}

We are now in a position to prove the main result of this paper,
Theorem~2 in Section~\ref{sec-const}.  I will address the five claims
in the theorem in order.  The proof of the final and principal claim
(e), that the modified chain has asymptotic variance at least as small
as the original chain, will follow closely the proof of Peskun's
theorem presented in the previous section.

\subsection*{Claim (a):\ \ The modified chain, with transition
  probabilities $\ddot T'$ is irreducible}\vspace*{-7pt}

Let $(a,b)$ and $(c,d)$ be distinct states in $\ddot \X$.  We need to 
show that $(a,b)$ and $(c,d)$ are linked by transitions with non-zero
probability under $\ddot T'$.
From the definition of $\ddot\X$,
$T(c,d)>0$.  If $b=c$, then $\ddot T'((a,b),(c,d))\, =\, U'_c(a,d)\, \ge\,
T(c,d)\, >\, 0$.  Otherwise, from the irreducibility of $T$, there exist
states $x_1,\ldots,x_k$ in $\X$ with $T(b,x_1)>0$, $T(x_k,c)>0$, and
$T(x_i,x_{i+1})>0$ for $i=1,\ldots,k\!-\!1$, and hence
$(b,x_1),\,(x_1,x_2),\,\ldots,\,(x_k,c) \,\in\, \ddot \X$.  Furthermore,
\beq
&&\!\!\!\!\!\!
   \ddot T'((a,b),(b,x_1))\ =\ U'_b(a,x_1)\ \ge\ T(b,x_1)\ >\ 0 \\[4pt]
&&\!\!\!\!\!\!
   \ddot T'((b,x_1),(x_1,x_2))\ =\ U'_{x_1}(b,x_2)\ \ge\ T(x_1,x_2)\ >\ 0 \\
&&\!\!\!\!\!\! 
   \mbox{\hspace*{1.2in}}\cdots\nonumber\\
&&\!\!\!\!\!\! 
   \ddot T'((x_k,c),(c,d))\ =\ U'_c(x_k,d)\ \ge\ T(c,d)\ >\ 0
\eeq

\subsection*{Claim (b):\ \ The transition probabilities 
$\ddot T'$ leave $\ddot \pi$ invariant}\vspace*{-7pt}

This is implied by the way $\ddot T'$ was constructed in
Section~\ref{sec-const}.  It can also be shown directly as follows:
\beq
  \sum_{(x_0,y_0)\in\ddot\X}\!\!\!
    \ddot\pi(x_0,y_0)\, \ddot T'((x_0,y_0),(x_1,y_1))
  & = &  \!\!\sum_{(x_0,y_0)\in\ddot\X}\!\!\!
    \pi(x_0)\,T(x_0,y_0)\,\delta(y_0,x_1)\,U'_{x_1}(x_0,y_1)\ \ \\[4pt]
  & = & \sum_{x_0\in\X} \pi(x_0)\, T(x_0,x_1)\,U'_{x_1}(x_0,y_1) \\[4pt]
  & = & \sum_{x_0\in\X} \pi(x_1)\, T(x_1,x_0)\, U'_{x_1}(x_0,y_1) \\[4pt]
  & = & \sum_{x_0\in\X} \pi(x_1)\, T(x_1,y_1)\, U'_{x_1}(y_1,x_0) \\[4pt]
  & = & \pi(x_1)\, T(x_1,y_1)\,\sum_{x_0\in\X} U'_{x_1}(y_1,x_0) \\[4pt]
  & = & \pi(x_1)\, T(x_1,y_1) \ \ =\ \ \ddot\pi(x_1,y_1)
\eeq

\subsection*{Claim (c):\ \ If $\X$ contains at least three elements, 
   $\ddot T'$ is not reversible}\vspace*{-7pt}

Let $a$, $b$, and $c$ be three distinct elements of $\X$.  Since the
original chain is irreducible, either $T(a,b)>0$ or there exist distinct
$x_1,\ldots,x_n$ such that $T(a,x_1)>0$, $T(x_n,b)>0$ and $T(x_i,x_{i+1})>0$
for $i=1,\ldots,n\!-\!1$.  Similarly, either $T(b,c)>0$ or there exist
distinct $y_1,\ldots,y_m$ linking $b$ to $c$.  One way or another, we
can find distinct $x,y,z$ such that $T(x,y)>0$ and 
$T(y,z)>0$ --- if $T(a,b)=0$, take three consecutive states from 
$a,x_1,\ldots,x_n,c$; if $T(b,c)=0$, take three consecutive states
from $b,y_1,\ldots,y_m,c$; and if $T(a,b)>0$ and $T(b,c)>0$, 
use $a,b,c$.  The states $(x,y)$ and $(y,z)$ are in $\ddot\X$, and 
have positive probability under $\ddot\pi$.  (Note that all states in $\X$
have positive probability under $\pi$, since the original chain is
irreducible.)  It follows that $\ddot T'((x,y),(y,z))$ is positive.
However, $\ddot T'((y,z),(x,y))$ is zero, since $z \ne x$.  The
modified chain with transition probabilities $\ddot T'$ is therefore
non-reversible.

\subsection*{Claim (d):\ \ The bias of the estimator $\hat\mu'_n$ is of order 
   $1/n$}\vspace*{-7pt}

The modified chain leaves invariant the distribution
$\ddot\pi(x,y)\,=\,\pi(x)T(x,y)\,=\,\pi(y)T(y,x)$.  The marginal distribution 
for the second component of state under $\ddot\pi$ is 
$\sum_x \ddot\pi(x,y) \,=\, \pi(y)$.  An MCMC estimator that looks at
a function of the second component of state will therefore converge to the 
correct expectation of this function with respect to $\pi$, with bias
of order $1/n$, in accordance with the standard properties of MCMC 
estimators, as discussed in Section~\ref{sec-prelim}.

\subsection*{Claim (e):\ \ The asymptotic variance of $\hat\mu'$ is no greater 
  than that of $\hat\mu$}\vspace*{-7pt}

This is the principal claim.  Its proof will follow the same steps as
the proof of Peskun's theorem in Section~\ref{sec-peskun-proof}.  

\subsection*{Step 1:\ \ Looking at one pair of states is enough}\vspace*{-7pt}

In Section~\ref{sec-const}, a chain on $\ddot\X$ was defined that was
essentially equivalent to the original chain on $\X$, with transitions
$T$.  The transition probabilities for this chain were defined in
equation~(\ref{eq-ddotT}) to have the form $\ddot T((x_0,y_0),\,(x_1,y_1)) \ =\
\delta(x_1,y_0)\,T(x_1,y_1)$, which can be seen as an instance of the
definition of $\ddot T'$ in equation~(\ref{eq-ddotTprime}), with 
$U_{x_1}(x_0,y_1)\ =\ T(x_1,y_1)$.  We can therefore view Theorem~2
as claiming that changing from this $U$ to some other $U'$ that satisfies
conditions~(\ref{eq-Ucond1}) and~(\ref{eq-Ucond2}) will not increase asymptotic
variance.  More generally, we will see that a change from transitions
$\ddot T$ based on \textit{any} $U$ satisfying
condition~(\ref{eq-Ucond1}) to transitions $\ddot T'$ 
based on some other $U'$ that also satisfies this condition and for which 
\beq 
  U'_x(y,z) \ \ge\ U_x(y,z),\ \ \ \mbox{for all $x,y,z\in\X$ with $y \ne z$}
\eeq
will not increase asymptotic variance.

Any such change from $U$ to $U'$ can be expressed as a sequences of changes, 
each of which affects $U_x(y,z)$ only when $x$ is some particular state
$O$, and $y$ and $z$ are both either state $A$ or state $B$.  In order
for condition~(\ref{eq-Ucond1}) to be satisfied, such $U$ and $U'$ must
be related as follows:
\beq
  U'_x(y,z)\ =\ U_x(y,z),\ \ \ 
     \mbox{when $x \ne O$ or $x \notin \{A,B\}$ or $y \notin \{A,B\}$} 
\nonumber\\
\label{eq-Urels}\\[-10pt]
    \begin{array}{ll} 
        U'_O(A,A)\ =\ U_O(A,A)-\delta_A,\ \ & U'_O(A,B)\ =\ U_O(A,B)+\delta_A 
     \\[2pt]
        U'_O(B,A)\ =\ U_O(B,A)+\delta_B,\ \ & U'_O(B,B)\ =\ U_O(B,B)-\delta_B
    \end{array}
\nonumber
\eeq
The resulting transition probabilities, $\ddot T$ and $\ddot T'$, will be
identical except for the following states:
\beq
   \ddot T'((A,O),(O,A))  =  \ddot T((A,O),(O,A)) - \delta_A,\ \ \
   \ddot T'((A,O),(O,B))  =  \ddot T((A,O),(O,B)) + \delta_A\ \ \nonumber
   \\[-8pt] \\[-8pt]
   \ddot T'((B,O),(O,A))  =  \ddot T((B,O),(O,A)) + \delta_B,\ \ \
   \ddot T'((B,O),(O,B))  =  \ddot T((B,O),(O,B)) - \delta_B\ \ \nonumber
\label{eq-Tprimerels}
\eeq
The rest of the proof will assume that $\ddot T$ and $\ddot T'$ differ
only as above.  The chain with transition probabilities $\ddot T$ will
be referred to as the ``old'' chain, while that with transitions probabilities
$\ddot T'$ will be called the ``new'' chain.

\subsection*{Steps 2 \& 3:\ \ Defining blocks delimited by ``delta'' 
   transitions}\vspace*{-7pt}

We can mark the transitions in a realization of either the old or the
new the chain that would have been different in the other chain ---
ie, those transitions where the addition or subtraction or $\delta_A$
or $\delta_B$ in~(\ref{eq-Tprimerels}) would have made a difference.
This can be done using a simulation procedure entirely analogous to
that described in Step~2 of the proof of Peskun's theorem.

As in Step~3 of that proof, we can use these ``delta'' transitions to
define the boundaries between ``blocks'' of states in a realization of
the old or new chain.  Such blocks will always begin with state
$(O,A)$ or $(O,B)$ and end with state $(A,O)$ or $(B,O)$.  If we
assume that we start the chain in state $(O,A)$ or $(O,B)$, a typical
sequence of blocks for the old chain might look like\\[8pt]
  \centerline{\psfig{file=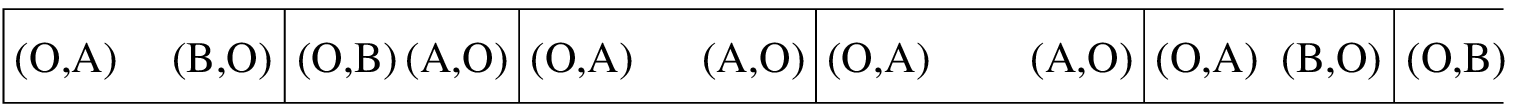,height=0.32in}}\\[8pt]
whereas a typical block sequence for the new chain might look like\\[8pt]
  \centerline{\psfig{file=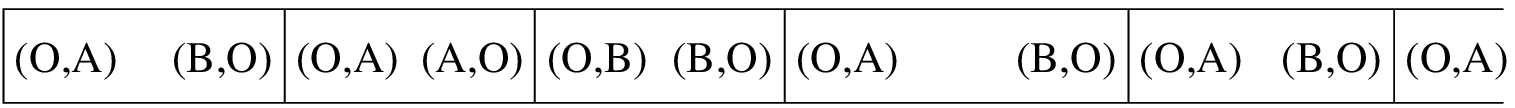,height=0.32in}}\\[8pt]
Note that states $(O,A)$, $(A,O)$, $(O,B)$, and $(B,O)$ may occur within
a block, as well as at the beginning and end.
The difference between the two chains is that states in the old chain on 
each side of a block boundary are simply reversals of each other, whereas 
in the new chain, one of the states on the two sides of a boundary 
will contain $A$ and the other $B$.

\subsection*{Step 4: 
  Probabilities of the four types of blocks and their contents}\vspace*{-7pt}

Blocks starting with $(O,A)$ and ending with $(A,O)$ will be called
$AA$ blocks, those starting with $(O,A)$ and ending with $(B,O)$ will
be called $AB$ blocks, and similarly for blocks of types $BB$ and $BA$.
The way these blocks are produced is illustrated in the diagram below:\\[12pt]
  \centerline{\psfig{file=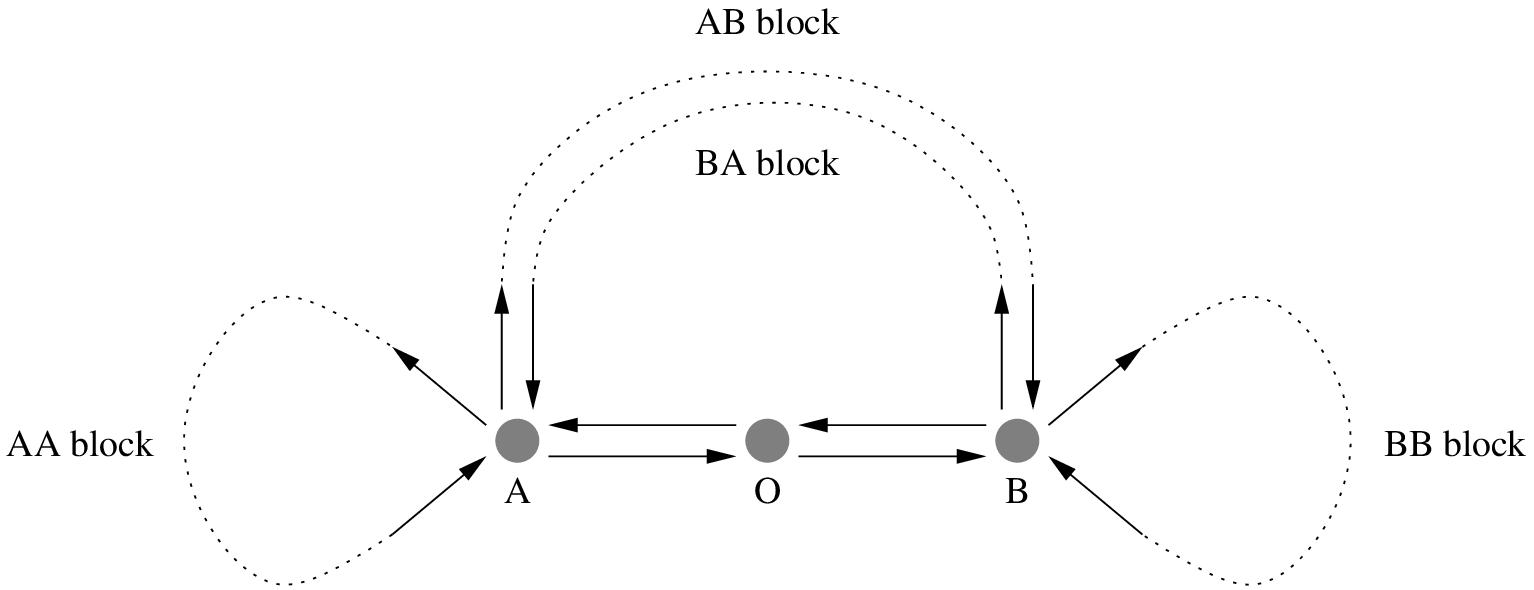}}\\[8pt]
An $AA$ block starts with state $(O,A)$, represented in the diagram 
by an arrow from $O$ to $A$.  Transitions from this state lead to an
arrow from $A$ to some other state, then arrows between other states,
and eventually to an arrow from some state to $A$, followed by an 
arrow from $A$ to $O$, which represents the the state $(A,O)$.  For the
block to end here, a delta transition must occur at this point, which
happens with probability $\delta_A$.  An $AB$ block also starts with an
arrow from $O$ to $A$, but ends with and arrow pointing to state $B$,
and then an arrow from $B$ to $O$ that is followed by a delta transition.

To find how the probabilities of the four types of blocks, $P(AA)$,
$P(AB)$, $P(BA)$, and $P(BB)$ are related, we can start by noting that
since condition~(\ref{eq-Ucond1}) applies to both $U$ and $U'$, which
are related by~(\ref{eq-Urels}), it follows that
\beq
    T(O,A)\,\delta_A & = & T(O,A)\ (U'_O(A,B) - U_O(A,B)) \nonumber\\[4pt]
                     & = & T(O,B)\ (U'_O(B,A) - U_O(B,A))
                     \ \ =\ \ T(O,B)\,\delta_B
\label{eq-rrel}
\eeq

Next, note that the probability that a block beginning with $(O,A)$ starts
at time $t$ (with $t$ being large) in the old chain is
\beq
  \ddot\pi(A,O))\,\delta_A & = & \pi(A)\,T(A,O)\,\delta_A \ \ =\ \ 
                                 \pi(O)\,T(O,A)\,\delta_A 
\eeq
Using (\ref{eq-rrel}), we see that this is equal to the probability that a 
block beginning with $(O,B)$ starts at time $t$, which is
\beq
  \ddot\pi(B,O)\,\delta_B & = & \pi(B)\,T(B,O)\,\delta_B \ \ =\ \
                                \pi(O)\,T(O,B)\,\delta_B
\eeq
Similarly, the probability that a block ends with state $(A,O)$ at time $t$,
\beq
  \ddot\pi(A,O)\,\delta_A & = & \pi(A)\,T(A,O)\,\delta_A \ \ =\ \ 
                                \pi(O)\,T(O,A)\,\delta_A
\eeq
is equal to the probability that a block ends with state $(B,O)$ at time $t$,
\beq
  \ddot\pi(B,O)\,\delta_B & = & \pi(B)\,T(B,O)\,\delta_B \ \ =\ \ 
                                \pi(O)\,T(O,B)\,\delta_B
\eeq
It follows that $P(AA)+P(AB)=P(BB)+P(BA)$ and $P(AA)+P(BA)=P(BB)+P(AB)$ 
for the old chain, from which we can conclude that $P(AA)=P(BB)$ and
$P(AB)=P(BA)$.  A similar argument shows this for the new chain as well.

The distribution of the contents of blocks of type $AA$ may differ
from that for the contents of blocks of type $BB$, but blocks of type
$AB$ and blocks of type $BA$ can be viewed as reversals of each other.
Consider, for example, a block consisting of the following states:
\beq
    (O,A),\ (A,X),\ (X,Y),\ (Y,B),\ (B,O)
\label{eq-exblock}
\eeq
As seen above, the probability of a block starting with $(O,A)$ at some
given time is $\pi(A)\,T(A,O)\,\delta_A$ (in both the old and the new chain,
as seen from~(\ref{eq-rrel}) and the reversibility of $T$).  
Multiplying by the probabilities
of the subsequent transitions, and the probability of a delta transition
from $(B,O)$, the probability of the $AB$ block above occurring at a given
time in the old chain is
\beq
  \pi(A)\,T(A,O)\,\delta_A\,U_A(O,X)\,U_X(A,Y)\,U_Y(X,B)\,U_B(Y,O)\,\delta_B
\label{eq-ABprob}
\eeq
This is also the probability in the new chain, since $U$ and $U'$ are the
same for all except delta transitions.
Now consider the reversal of the block in~(\ref{eq-exblock}):
\beq
    (O,B),\ (B,Y),\ (Y,X),\ (X,A),\ (A,O)
\label{eq-exblockrev}
\eeq
The probability of this block occurring at a given time is
\beq
   \pi(B)\,T(B,O)\,\delta_B\,U_B(O,Y)\,U_Y(B,X)\,U_X(Y,A)\,U_A(X,O)\,\delta_A
\label{eq-BAprob}
\eeq
We can see that (\ref{eq-ABprob}) and (\ref{eq-BAprob}) are equal as follows,
using the reversibility of $T$ with respect to $\pi$ and the fact that
$U$ satisfies condition~(\ref{eq-Ucond1}):
\beq
  \lefteqn{\pi(A)\,T(A,O)\,U_A(O,X)\,U_X(A,Y)\,U_Y(X,B)\,U_B(Y,O)
  }\ \ \ \ \
  \nonumber\\[4pt]
  & = &
  \pi(A)\,T(A,X)\,U_A(X,O)\ \cdot\ U_X(A,Y)\,U_Y(X,B)\,U_B(Y,O) \\[4pt]
  & = &
  U_A(X,O)\ \cdot\ \pi(X)\,T(X,A)\,U_X(A,Y)\ \cdot\ U_Y(X,B)\,U_B(Y,O) \\[4pt]
  & = &
  U_A(X,O)\ \cdot\ \pi(X)\,T(X,Y)\,U_X(Y,A)\ \cdot\ U_Y(X,B)\,U_B(Y,O) \\[4pt]
  & = &
  U_X(Y,A)\,U_A(X,O)\ \cdot\ \pi(Y)\,T(Y,X)\,U_Y(X,B)\ \cdot\ U_B(Y,O) \\[4pt]
  & = &
  U_X(Y,A)\,U_A(X,O)\ \cdot\ \pi(Y)\,T(Y,B)\,U_Y(B,X)\ \cdot\ U_B(Y,O) \\[4pt]
  & = &
  U_Y(B,X)\,U_X(Y,A)\,U_A(X,O)\ \cdot\ \pi(B)\,T(B,Y)\,U_B(Y,O) \\[4pt]
  & = &
  U_Y(B,X)\,U_X(Y,A)\,U_A(X,O)\ \cdot\ \pi(B)\,T(B,O)\,U_B(O,Y) \\[4pt]
  & = &
  \pi(B)\,T(B,O)\,U_B(O,Y)\,U_Y(B,X)\,U_X(Y,A)\,U_A(X,O)
\eeq

\subsection*{Step 5: Sampling in the new chain is stratified}\vspace*{-7pt}

As was done in the proof of Peskun's theorem, we can now imagine
simulating both the old chain and the new chain one block at a time.
For each block, we decide whether it should be homogeneous (of type
$AA$ or $BB$) or non-homogeneous (of type $AB$ or $BA$).  The
probability of a block being homogeneous is the same regardless of
whether the block starts with $A$ or $B$, since the results in Step~4
imply that $P(AA)\,/\,(P(AA)+P(AB))\, = \,P(BB)\,/\,(P(BB)+P(BA))$.  The
probability that a block is homogeneous is also the same for the old
chain and the new chain.  If we make the same random decisions as to
whether or not blocks are homogeneous in the old and new chains, the
sequence of homogeneous versus non-homogeneous blocks will be the same
for the two chains.

The only significant difference between the old and new chains is that
in the new chain the sequence of homogeneous blocks alternates between
$AA$ blocks and $BB$ blocks, and hence the number of $AA$ blocks is
equal to the number of $BB$ blocks (or differs by only one).  This
arises for exactly the same reasons as in the proof of Peskun's
theorem --- if $N_{AA}$, $N_{AB}$, $N_{BA}$, and $N_{BB}$ are the numbers 
of blocks of each type, $|\, (N_{AA}\!+\!N_{BA})\, -\, (N_{BB}+N_{BA})\,|
 \,\le\, 1$ in the new chain, and hence $|N_{AA}-N_{BB}| \le 1$.
The sampling for $AA$ and $BB$ blocks is therefore
stratified in the new chain, but not in the old chain.  The old chain
stratifies sampling of $AB$ and $BA$ blocks, but since the
distributions for the contents of blocks of types $AB$ and $BA$ are
the same (apart from a reversal, which doesn't affect sums of function
values), this stratification in the old chain has no effect.

\subsection*{Steps 6: Stratification will not increase asymptotic 
variance}\vspace*{-7pt}

Finally, as in the proof of Peskun's theorem, we can apply Lemma~1 in
the Appendix to show that the asymptotic variance using a simulation
that continues for a specified number of blocks is the same as that
using a simulation for a specified number of transitions.  We can then
apply Lemma~2 to show that the block-by-block simulation of the new
chain, which is stratified with respect to $AA$ and $BB$ blocks, will
have asymptotic variance at least as small as for the old chain.

\section{\hspace*{-7pt}Conclusion}\label{sec-conc}\vspace*{-10pt}

This paper shows how any reversible Markov chain can be transformed
into a non-reversible chain that tries to avoid backtracking to the
state visited immediately before.  This transformation never increases
the asymptotic variance of an MCMC estimator using the chain, and will
usually decrease it.  Sometimes, the decrease in asymptotic variance
is dramatic, but other times it is small.  In general, one would
expect the decrease in asymptotic variance to be small when a state of
the original Markov chain has many possible successor states (of
roughly similar probability), since in this situation, even the
original chain will rarely backtrack.  In many circumstances, the
chain that avoids backtracking will require little or no more time per
transition than the original chain, though this cannot be guaranteed
in all cases.

The particular transformation described in this paper may sometimes be
of practical use.  For many problems, however, the gains may be slight
or non-existent.  In particular, for problems with continuous state
spaces, and continuous transition distributions, exact backtracking
has zero probability of occurring anyway.  There seems to be scope for
generalizing the idea of trying to avoid backtracking, however.
Possibilities include trying to avoid backtracking to any of the past
several states, and trying to avoid backtracking not just to the exact
previous state, but also to anywhere in its vicinity.

More generally, the results in this paper indicate that non-reversible
Markov chains have a fundamental advantage over reversible chains, and
that the search for better MCMC methods may therefore be best focused
on non-reversible chains.  The proof techniques used in this paper may
be useful in analysing such methods.

\section*{Acknowledgements}\vspace{-10pt}

I thank Longhai Li and Jeffrey Rosenthal for helpful discussions.
This research was supported by the Natural Sciences and Engineering
Research Council of Canada.  The author holds a Canada Research Chair
in Statistics and Machine Learning.

\section*{Appendix:  Statements and proofs of lemmas}\vspace*{-10pt}

The following lemmas were used in the proofs of Sections~\ref{sec-peskun-proof}
and~\ref{sec-main-proof}. 

The first lemma justifies looking at the asymptotic variance of
simulations continuing for a specified number of blocks, instead of a
specified number of transitions.  To apply it to blocks defined by
``delta'' transitions, we can extend the state space of the Markov
chain to include an indicator of whether the current state is the last
in a block (essentially moving the decision whether a transition from
state $A$ or $B$ is to be ``marked'' back to the previous transition
into state $A$ or $B$).  With this extension, the set $\S$ below can
consist of states $A$ or $B$ at the end of a block.

\vspace*{6pt}

{\em

\noindent\textbf{Lemma 1:}\ \ Let $X_1,X_2,\ldots$ be an irreducible
Markov chain on a finite state space $\X$, with invariant
distribution $\pi(x)$. Let $\S$ be some non-empty subset of $\X$, and let
$f(x)$ be some function of state, whose expectation with respect to 
$\pi$ is $\mu$.
Define\vspace*{-8pt}
\beq 
 \textstyle N(k) & = & \min\, \Big\{n : \sum\limits_{t=1}^n I_{\S}(X_t)\, =\, k
                       \Big\}
\eeq
where $I_{\S}$ is the indicator function for $\S$.  Consider the following 
two families of estimators:\vspace*{-5pt}
\beq
   \hat\mu_n = {\displaystyle {1 \over n}} \sum\limits_{t=1}^n f(X_t),\ \ \ \ \
   \tilde\mu_k = {\displaystyle {1 \over N(k)}} \sum\limits_{t=1}^{N(k)} f(X_t)
\eeq
The asymptotic variances of these estimators are the same:\vspace*{-2pt}
\beq 
   \lim_{n\rightarrow\infty} n\, \Var(\hat\mu_n) \ =\ 
   \lim_{n\rightarrow\infty} \, n\, \Var(\tilde\mu_{\lceil n\pi(\S) \rceil})
\eeq

}

\noindent \textbf{Proof:}  Without loss of generality, suppose $\mu=0$.
We will see that as $n$ increases, $n\Var(\hat\mu_n)$ and 
$n\Var(\tilde\mu_{\lceil n\pi(\S) \rceil})$ both approach $(n+n^{1/2+\epsilon})
\Var(\hat\mu_{n+n^{1/2+\epsilon}})$, where $\epsilon$ is a positive
constant to be set below.  The diagram below may help to visualize the 
proof:\vspace{1pt}

\begin{minipage}{6.5in}

\hspace*{1.00in} $0$ \hspace*{1.85in} $n-n^{1/2+\epsilon}$ \hspace*{0.32in} $n$
                 \hspace*{0.45in} $n+n^{1/2+\epsilon}$ \\
\hspace*{1.03in} \rule{4in}{1pt} \\[-10pt]
\hspace*{1.03in} \rule{1pt}{10pt} \hspace*{2.15in} \rule{1pt}{10pt} 
                 \hspace*{0.75in} \rule{1pt}{10pt} 
                 \hspace*{0.75in} \rule{1pt}{10pt} \\[-12pt]
\hspace*{3.85in} \rule{1pt}{20pt} \\
\hspace*{3.45in} $N(\lceil n\pi(\S) \rceil)$

\end{minipage}

\vspace{1pt}

First, we note that $(n+n^{1/2+\epsilon})\,\hat\mu_{n+n^{1/2+\epsilon}} =
n \hat\mu_n + n^{1/2+\epsilon} Z$, where $Z$ is the average 
of $f(X_i)$ for $i$ from $n+1$ to $n+n^{1/2+\epsilon}$.  Dividing
by $\sqrt{n+n^{1/2+\epsilon}}$, we get
\beq 
   \sqrt{n+n^{1/2+\epsilon}}\,\hat\mu_{n+n^{1/2+\epsilon}} =
   \sqrt{n/(n+n^{1/2+\epsilon})} \Big[ \sqrt{n} \hat\mu_n + n^{\epsilon} Z\Big]
\eeq
As $n$ increases, the first factor on the right will go to one.  By
the Central Limit Theorem for Markov chains, $|Z|$ will be less than
$(n^{1/2+\epsilon})^{-1/2+\epsilon} = n^{-1/4+\epsilon^2}$ with
probability approaching one exponentially fast, so if $\epsilon$ is in
$(0,\,(\sqrt{2}\!-\!1)/2)$, the term $n^{\epsilon}Z$ will go to zero.
It follows that $n\Var(\hat\mu_n)$ will approach
$(n+n^{1/2+\epsilon})\Var(\hat\mu_{n+n^{1/2+\epsilon}})$.  (Since
$f(x)$ is bounded, an exponentially small probability of a large value
for $|Z|$ cannot affect this limit.)

From the Central Limit Theorem, we can also conclude that 
$N(\lceil n\pi(\S) \rceil)$ will be in the
interval $(n\!-\!n^{1/2+\epsilon},\,n\!+\!n^{1/2+\epsilon})$ with probability
approaching one exponentially fast.  If so, we can write 
\beq
   (n\!+\!n^{1/2+\epsilon})\,\hat\mu_{n+n^{1/2+\epsilon}}\ =\
   N(\lceil n\pi(\S) \rceil)\, \tilde\mu_{\lceil n\pi(\S) \rceil}
     \,+\, (n\!+\!n^{1/2+\epsilon}\!-\!N(\lceil n\pi(\S) \rceil))\, Y
\eeq
where $Y$ is the average 
of $f(X_i)$ for $i$ from $N(\lceil n\pi(\S) \rceil)+1$ to $n+n^{1/2+\epsilon}$.
Dividing by $\sqrt{n+n^{1/2+\epsilon}}$, we get\vspace*{-6pt}
\beq
   \sqrt{n+n^{1/2+\epsilon}}\,\hat\mu_{n+n^{1/2+\epsilon}} =
   {N(\lceil n\pi(\S) \rceil) \over n \sqrt{1+n^{-1/2+\epsilon}}}
   \Big[ \sqrt{n} \tilde\mu_{\lceil n\pi(\S) \rceil}
    + (\sqrt{n} / N (\lceil n\pi(\S) \rceil)) \,KY \Big]
\eeq
where $K = n\!+\!n^{1/2+\epsilon}\!-\!N(\lceil n\pi(\S) \rceil)$ will
be in $(0,2n^{1/2+\epsilon})$ if $N(\lceil n\pi(\S) \rceil)$ is in
$(n\!-\!n^{1/2+\epsilon},\,n\!+\!n^{1/2+\epsilon})$.\pagebreak\linebreak
By the Central Limit Theorem for
Markov chains, $|KY|$ will be less than $(2n^{1/2+\epsilon})^{1/2+\epsilon}
= 2^{1/2+\epsilon}n^{1/4+\epsilon+\epsilon^2}$ with a probability that
approaches one exponentially fast.  Since $N(\lceil n\pi(\S) \rceil)$ will 
approach $n$, we can see that $n\Var(\tilde\mu_{\lceil n\pi(\S) \rceil})$
will approach $(n+n^{1/2+\epsilon})\Var(\hat\mu_{n\!+\!n^{1/2+\epsilon}})$.

Since $n\Var(\hat\mu_n)$ and $n\Var(\tilde\mu_{\lceil n\pi(\S)
\rceil})$ both approach the same value as $n$ goes to infinity, they
must also have the same limit, as the lemma states.

\vspace*{6pt}

The second lemma justifies the claim that partial stratification of
sampling for blocks cannot increase asymptotic variance.  In applying
this lemma to the proofs in Sections~\ref{sec-peskun-proof}
and~\ref{sec-main-proof}, $Z_1,Z_2,\ldots$ are identifiers for the
type of each block (we can use $0=AA$, $1=BB$, and
\mbox{$2=AB\mbox{~or~}BA$}), which form a Markov chain, since the
distribution for the type of the next block depends only on the type
of the previous block.  The types of blocks for the modified chain are
$Z'_1,Z'_2,\ldots$, which are stratified with respect to $0$ and $1$.
In the applications of this lemma, $H$ corresponds to the sum of the
values of $f$ for all states in a block, and $L$ corresponds to the
number of states in this block.

\vspace*{6pt}

{\em 

\noindent\textbf{Lemma 2:}\ \ Let $Z_1,Z_2,\ldots$ be an irreducible
Markov chain with state space $\{ 0, 1, 2 \}$, whose invariant
distribution, $\rho$, satisfies $\rho(0)=\rho(1)$.  Let $Q_z$ for 
$z=0,1,2$ be distributions for
pairs $(H,L) \in {\mathbb R} \times {\mathbb R}^+$ having
finite second moments.  Conditional on
$Z_1,Z_2,\ldots$, let $(H_i,L_i)$ be drawn independently from $Q_{Z_i}$.  
Define\vspace*{-10pt}
\beq
   Z'_i & = & \left\{\begin{array}{ll} 
    Z_i & \mbox{if $Z_i=2$} \\[1pt]
    Z_k + \sum\limits_{j=1}^{i-1} I_{\{0,1\}}(Z_j)\ \ (\mbox{\rm modulo 2}) 
      & \mbox{if $Z_i\ne2$}
   \end{array}\right.
\eeq
where $k=\min\{i:Z_i\ne2\}$.  (In other words, the $Z'_i$ are the same
as the $Z_i$ except that the positions where $0$ or $1$ occurs have
their values changed to a sequence of alternating $0$s and $1$s.)
Conditional on
$Z_1,Z_2,\ldots$, let $(H'_i,L'_i)$ be drawn independently from $Q_{Z'_i}$.  
Define two families of estimators as follows:\vspace*{-10pt}
\beq
   R_n \,=\, \sum_{i=1}^n H_i \,\Big/\, \sum_{i=1}^n L_i,\ \ \ \ \
   R'_n \,=\, \sum_{i=1}^n H'_i \,\Big/\, \sum_{i=1}^n L'_i
\eeq
Then the asymptotic variance of $R'$ is no greater than that of $R$.  In other
words, 
\beq 
   \lim\limits_{n\rightarrow\infty} n\Var(R'_n) & \le &
   \lim\limits_{n\rightarrow\infty} n\Var(R_n)
\eeq

}

\noindent\textbf{Proof:}
Let $N_{n,m} = (1/n)\sum_{i=1}^n I_{\{m\}}(Z_i)$ and $N'_{n,m} = (1/n)
\sum_{i=1}^n I_{\{m\}}(Z'_i)$.  Note that $E(N_{n,m})=E(N'_{n,m})$ and
$|N'_{n,1}\!-\!N'_{n,0}| \le 1/n$, so the proportions of pairs from $Q_0$
and $Q_1$ are stratified in $R'_n$. 
By the Central Limit Theorem for Markov chains, $N_n=(N_{n,0},N_{n,1},N_{n,2})$
and $N'_n = (N'_{n,0},N'_{n,1},N'_{n,2})$ asymptotically have (degenerate)
multivariate normal distributions, with the same mean vectors, though
different covariance matrices.  (Note that although $Z'_1,Z'_2,\ldots$ is
not a Markov chain, a Markov chain can be defined on an extended state space 
that includes $Z'_i$ as a component.)

The asymptotoic variances of $R_n$ and $R'_n$ can be decomposed as follows:
\beq
  n\Var(R_n) \ =\ n\Var(E(R_n|N_n)) \ +\ E(n\Var(R_n|N_n)) \label{eq-V1}\\[2pt]
  n\Var(R'_n) \ =\ n\Var(E(R'_n|N'_n)) \ +\ E(n\Var(R'_n|N'_n))\label{eq-V2}
\eeq
Note that the distribution of $R_n$ given $N_n=N$ is the same as the
distribution of $R'_n$ given $N'_n=N$.  Writing $R_n = (1/n)\sum H_i\ /\ (1/n)
\sum L_i$, we can apply the Central Limit Theorem to the numerator and
denominator, then use the delta rule to conclude that $R_n$ given $N_n$
is asymptotically normal, with asymptotic variance that depends only on $N_n$
(not on $n$).
Since $N_n$ and $N'_n$ have the same means and both are asymptotically normal, 
we can apply the delta rule again to conclude that the second term on the 
right in~(\ref{eq-V1}) is equal to the second term on the right 
in~(\ref{eq-V2}).

Looking at the first terms in~(\ref{eq-V1}) and~(\ref{eq-V2}), we can
rewrite $n\Var(E(R_n|N_n))$ and $n\Var(E(R'_n|N'_n))$ as follows:\vspace*{-8pt}
\beq
 n\Var(E(R_n|N_n)) & = & n\Var(E(E(R_n|N_n)|N_{n,2}))
                          \ +\ E(n\Var(E(R_n|N_n)|N_{n,2})) 
\label{eq-VV1} \\[4pt]
 n\Var(E(R'_n|N'_n)) & = & n\Var(E(E(R'_n|N'_n)|N'_{n,2})) 
                          \ +\ E(n\Var(E(R'_n|N'_n)|N'_{n,2}))
\label{eq-VV2}
\eeq
Since the expectations of $N_{n,1}$ and $N_{n,2}$ given $N_{n,2}=N$
are the same as the expectations of $N'_{n,1}$ and $N'_{n,2}$ given 
$N'_{n,2}=N$, we can conclude that $E(E(R_n|N_n)|N_{n,2}=N)$ is asymptotically
equal to $E(E(R'_n|N'_n)|N'_{n,2}=N)$.  The distributions of $N_{n,2}$
and $N'_{n,2}$ are the same, so it follows that the first terms on the
right in~(\ref{eq-VV1}) and~(\ref{eq-VV2}) are asymptotically equal.
Due to stratification, $N'_{n,0}$ and $N'_{n,1}$ are fixed given 
$N'_{n,2}$, so that $\Var(E(R'_n|N'_n)|N'_{n,2})=0$.  It follows that 
the second terms on the right in~(\ref{eq-VV1}) and~(\ref{eq-VV2})
are related by
\beq
 E(n\Var(E(R_n|N_n)|N_{n,2})) & \ge & E(n\Var(E(R'_n|N'_n)|N'_{n,2})) \ \ =\ \ 0
\eeq

Combining these results, we can conclude that
asymptotically $\Var(E(R_n|N_n)) \ge \Var(E(R'_n|N'_n))$, and finally, that
$\Var(R_n)$ is asymptotically at least as large as $\Var(R'_n)$.

\section*{References}\vspace*{-10pt}

\leftmargini 0.2in
\labelsep 0in

\begin{description}
\itemsep 0pt

\item[]
  Adler, S.~L.\ (1981) ``Over-relaxation method for the Monte Carlo evaluation
  of the partition function for multiquadratic actions'', {\em Physical
  Review D}, vol.~23, pp.~2901-2904.

\item[]
  Diaconis, P., Holmes, S., and Neal, R.~M.\ (2000) ``Analysis of a 
  non-reversible Markov chain sampler'', {\em Annals of Applied Probability},
  vol.~10, pp.~726-752.

\item[]
  Gustafson, P.\ (1998) ``A guided walk Metropolis algorithm'', 
  \textit{Statistics and Computing}, vol.~8, pp.~357-364.

\item[]
  Hastings, W.~K.\ (1970) ``Monte Carlo sampling methods using Markov chains 
  and their applications'', {\em Biometrika}, vol.~57, pp.~97-109.

\item[]
  Hoel, P.~G., Port, S.~C., and Stone, C.~J.\ (1972) \textit{Introduction
  to Stochastic Processes}, Waveland Press.

\item
  Horowitz, A.~M.\ (1991) ``A generalized guided Monte Carlo algorithm'', 
  \textit{Physics Letters B}, vol.~268, pp.~247-252.

\item
  Kemeny, J.~G.\ and Snell, J.~L.\ (1960) \textit{Finite Markov Chains},
  Van Nostrand.

\item[]
  Peskun, P.~H.\ (1973) ``Optimum Monte-Carlo sampling using Markov chains'',
  \textit{Biometrika}, vol.~60, pp.~607-612.

\item[]
  Liu, J.~S.\ (1996) ``Peskun's theorem and a modified discrete-state
  Gibbs sampler'', \textit{Biometrika}, vol.~83, pp.~681-682.

\item[]
  Liu, J.~S.\ (2001) \textit{Monte Carlo Strategies in Scientific Computing},
  Springer-Verlag.

\item[]
  Mira, A.\ and Geyer, C.~J.\ (2000) ``On non-reversible Markov chains'',
  in N.~Madras (editor) \textit{Monte Carlo Methods}, Fields Institute / AMS,
  pp.~95-110.

\item[]
  Neal, R.~M.\ (1998) ``Suppressing random walks in
  Markov chain Monte Carlo using ordered overrelaxation'', in
  M.~I.~Jordan (editor) {\em Learning in Graphical Models}, 
  Dordrecht: Kluwer Academic Publishers.

\item[]
  Neal, R.~M.\ (2003) ``Slice sampling'' (with discussion), \textit{Annals
  of Statistics}, vol.~31, pp.~705-767.

\item[]
  Romanovsky, V.~I.\ (1970) \textit{Discrete Markov Chains}, translated
  from the Russian by E.~Seneta, Wolters-Noordhoff.

\item[]
  Tierney, L.\ (1994) ``Markov chains for exploring posterior distributions''
  (with discussion), \textit{Annals of Statistics}, vol.~22, pp.~1701-1762.

\item[]
  Tierney, L.\ (1998) ``A Note on Metropolis-Hastings kernels for general
  state spaces'', \textit{Annals of Applied Probability}, vol.~8, pp.~1-9.

\end{description}

\end{document}